\documentclass[pdflatex,sn-mathphys-num]{sn-jnl}

\usepackage{graphicx}%
\usepackage{multirow}%
\usepackage{amsmath,amssymb,amsfonts}%
\usepackage{amsthm}%
\usepackage{mathrsfs}%
\usepackage[title]{appendix}%
\usepackage{xcolor}%
\usepackage{textcomp}%
\usepackage{manyfoot}%
\usepackage{booktabs}%
\usepackage{algorithm}%
\usepackage{algorithmicx}%
\usepackage{algpseudocode}%
\usepackage{listings}%

\usepackage{array, siunitx}
\usepackage{float}

\numberwithin{equation}{section}
\numberwithin{table}{section}

\theoremstyle{thmstyleone}%
\newtheorem{theorem}{Theorem}[section]
\newtheorem{proposition}[theorem]{Proposition}%
\newtheorem{lemma}[theorem]{Lemma}
\newtheorem{remark}[theorem]{Remark}

\begin{document}


\title[Linear stability from local bifurcation structure]{Linear stability of the first bifurcation in a tumor growth free boundary problem via local bifurcation structure}


\author[1]{\fnm{Junying} \sur{Chen}}\email{chenjy685@mail2.sysu.edu.cn}

\author*[1]{\fnm{Ruixiang} \sur{Xing}}\email{xingrx@mail.sysu.edu.cn}

\affil[1]{\orgdiv{School of Mathematics}, \orgname{Sun Yat-sen University}, \orgaddress{ \city{Guangzhou}, \postcode{510275}, \country{P.R. China}}}


\abstract{In this paper, we consider a 3-dimensional free boundary problem modeling tumor growth with the Robin boundary condition. The system involves a positive parameter $\mu$ which reflects the intensity of tumor aggressiveness. 
Huang, Zhang and Hu [Nonlinear Anal. Real World Appl. 2017(35), 483-502] have shown that for each $\mu_n$ ($n$ even) in a strictly increasing sequence $\{ \mu_n \}(n\geq 2)$, there exists a stationary bifurcation solution $(\sigma_n(\varepsilon),p_n(\varepsilon),r_n(\varepsilon))$ with $\mu = \mu_n(\varepsilon)$ bifurcating from $\mu_n$. We first derive that the bifurcation curve $(r_2(\varepsilon),\mu_2(\varepsilon))$ exhibits a transcritical bifurcation with $\mu_2'(0)<0$. Moreover, we show that the stationary bifurcation solution $(\sigma_2(\varepsilon),p_2(\varepsilon),r_2(\varepsilon))$ is linearly unstable for small $|\varepsilon|$ under non-radially symmetric perturbations. In contrast to the linear stability of the radially symmetric stationary solution, the lack of explicit expressions for bifurcation solutions adds great difficulty in analyzing their linear stability. 
The novelty of this paper lies in the use of the bifurcation curve's structure to overcome the above difficulties and avoid a large number of calculations.}




\maketitle

\section{Introduction}\label{sec1}

Consider the following free boundary problem
\begin{align}
\label{1.1}
&c\frac{\partial \sigma}{\partial t}-\Delta \sigma+  \sigma=0, \quad && x \in \Omega(t), t>0,\\
\label{1.2}  
&\frac{\partial \sigma}{\partial \vec{n}}+\beta(t)(\sigma-1)=0, && x \in  \partial \Omega(t), t>0,  \\
\label{1.3} 
&-\Delta p=\mu(\sigma-\tilde{\sigma}), \quad && x \in \Omega(t), t>0,\\
\label{1.4} 
&p=\gamma \kappa, && x \in \partial \Omega(t), t>0,\\
\label{1.5} 
&V_n=-\frac{\partial p}{\partial \vec{n}}, \quad && x \in \partial \Omega(t), t>0, \\
\label{1.6}
&\partial \Omega(0) = \partial \Omega_0, \qquad \left.\sigma\right|_{t=0}  =\sigma_0,  
\end{align}
where $\Omega(t) \subset \mathbb{R}^3$ is the tumor domain. $\sigma$ and $p$ denote the nutrient concentration and the internal pressure within the tumor, respectively. $c$ is a small positive constant representing the ratio between the cell proliferation rate and the nutrient diffusion rate. $\vec{n}$ is the unit outward normal direction, $\beta(t)$ is a positive function representing the rate of nutrient supply and  $1$ is the simplified external nutrient concentration. $\mu$ is a positive constant reflecting the intensity of tumor aggressiveness and $\tilde{\sigma}>0$ is the threshold concentration. $\gamma$ denotes the cell-to-cell adhesiveness and $\kappa$ represents the mean curvature of $\partial \Omega(t)$. Let $V$ be the velocity of movement of tumor cells. Supposing that the tumor domain is a porous medium, Darcy's law $V=-\nabla p$ holds. Together with the conversation of mass $\operatorname{div} V = \mu (\sigma-\tilde{\sigma})$, we obtain \eqref{1.3}. Assuming that the velocity $V$ is continuous up to the boundary, \eqref{1.5} holds. The boundary condition \eqref{1.2} was first proposed by Friedman and Lam in \cite{FriedmanLam2015radial}.
They investigated the properties of the radius of spherical tumors and the asymptotic behavior of radially symmetric solutions under various assumptions on $\beta(t)$. Also, they established the existence and uniqueness of radially symmetric stationary solutions.

The corresponding stationary problem of \eqref{1.1}--\eqref{1.5} is
\begin{align}
\label{sigma equation}
&-\Delta \sigma+\sigma=0, \quad && x \in \Omega,\\
\label{sigma boundary} &\frac{\partial \sigma}{\partial  \vec{n}}+\beta(\sigma-1)=0, && x \in  \partial \Omega,\\
\label{p equation} 
&-\Delta p=\mu(\sigma-\tilde{\sigma}), \quad && x \in \Omega,\\
\label{p boundary} 
&p= \gamma \kappa, && x \in \partial \Omega,\\
\label{dn p=0} 
&\frac{\partial p}{\partial \vec{n}}=0, \quad && x \in \partial \Omega. 
\end{align}

In 2017, Huang, Zhang and Hu \cite{HuangZhangHu2017bifurcation} considered problem \eqref{sigma equation}--\eqref{dn p=0} and proved that there exists a strictly increasing sequence $\{ \mu_n \}(n\geq 2)$ such that $\mu_n$($n$ even) is a bifurcation value, in which the axial symmetry solution $(\sigma_n(\varepsilon),p_n(\varepsilon),r_n(\varepsilon))$ of \eqref{sigma equation}--\eqref{dn p=0} with $\mu=\mu_n(\varepsilon)$ bifurcating from the radially symmetric solution $(\sigma_s,p_s,R)$(given in \cite{FriedmanLam2015radial}) with the free boundary $r_n(\varepsilon) = R + \varepsilon Y_{n, 0}+ O(\varepsilon^2)$. Subsequently, for \eqref{1.1}--\eqref{1.6} with $c=1$ and $\beta(t)=\beta$, they studied the linear and asymptotic stability of the radially symmetric stationary solution under non-radially symmetric perturbations in \cite{HuangZhangHu2019linear,HuangZhangHu2021asymptotic}. The study of extended models involving the Robin boundary condition in free boundary problems of tumor growth is referred to \cite{ZhuangCui2018asymptotic,Xu2020time,SongWangHu2024necrotic,ZhangZhang2021periodic,PengFengWei2024inhibitor,WangSongXu2018inhibitor} and references therein. These results mainly focus on the properties and asymptotic behavior of radially symmetric solutions, the existence of radially symmetric stationary solutions and the existence of stationary bifurcation solutions, as well as the linear and asymptotic stability of radially symmetric stationary solutions under non-radially symmetric perturbations. However, they do not address the bifurcation structure and the linear stability of stationary bifurcation solutions. To date, the only result on this topic in the free boundary problems of tumor growth was presented by Friedman and Hu in \cite{FriedmanHu2008}. They studied problem \eqref{1.1}--\eqref{1.6} with $c=1$ and \eqref{1.2} replaced by $\sigma = 1$ and showed that there exists an increasing sequence $\{ \hat{\mu}_n \}$ ($n \geq$ $2$) such that a branch of axial symmetry stationary solutions $( \hat{\sigma}_n(\varepsilon), \hat{p}_n(\varepsilon),\hat{r}_n(\varepsilon))$ with $\mu=\hat{\mu}_n(\varepsilon)$ bifurcating from the radially symmetric stationary solution (given in \cite{FriedmanReitich1999radial}) at $\mu = \hat{\mu}_n$($n$ even). Furthermore, they derived that the first bifurcation curve $\{(\hat{r}_2(\varepsilon), \hat{\mu}_2(\varepsilon))\}$ is a transcritical bifurcation and 
\begin{align}
\label{cite dmu(0)}
\hat{\mu}_2'(0)= 
\frac{1}{28} \sqrt{\frac{5}{\pi}} \hat{\mu}_2^2  \hat{R}^4 P_0(\hat{R})\left(-\frac{1}{2} P_1(\hat{R})-\frac{1}{2} P_2(\hat{R})-\hat{R}^2 P_2^2(\hat{R})+\hat{R}^2 P_1(\hat{R}) P_2(\hat{R})\right)<0,
\end{align}
where 
\begin{align}
\label{P_n}
P_n(\hat{R})=\frac{I_{n+3/2}(\hat{R})}{\hat{R} I_{n+1/2}(\hat{R})},
\end{align}
and $I_m$ is the modified Bessel function.
Through a large number of skillful calculations and estimations, they obtained that for problem \eqref{1.1}--\eqref{1.5} with $c=1$ and \eqref{1.2} replaced by $\sigma=1$, the stationary bifurcation solution $(\hat{\sigma}_2(\varepsilon), \hat{p}_2(\varepsilon), \hat{r}_2(\varepsilon))$ is linearly stable for $\varepsilon>0$ and linearly unstable for $\varepsilon<0$ under non-radially symmetric perturbations.

For problem \eqref{1.1}--\eqref{1.5} with $c=0$ and $\beta(t)  = \beta$, we study the linear stability of the first bifurcation solution $(\sigma_2(\varepsilon),p_2(\varepsilon), r_2(\varepsilon))$ of problem \eqref{sigma equation}--\eqref{dn p=0} from the view of the local bifurcation structure. At first, we give the result about the local bifurcation structure. 

\begin{theorem}{\bf (Transcritical bifurcation)}
\label{result1}
Let $\{(r_2(\varepsilon), \mu_2(\varepsilon))\}$ be the bifurcation curve of problem \eqref{sigma equation}--\eqref{dn p=0}. Then 
\begin{align}
\label{dmu(0)}
\mu_2'(0) = \dfrac{1}{7} \sqrt{\dfrac{5}{\pi}} \frac{\mu_2^2 R^3}{4 \gamma} \frac{  \beta P_0(R) }{\beta+RP_{0}(R) }
\frac{1}{ \left( \beta+\dfrac{2}{R}+R P_2(R) \right)^2} \left( E_1 R\beta^2 + E_2 \beta + E_3 \right)<0,
\end{align}
where  
\begin{align}
\label{E_1}
E_1=&-\frac{1}{2}\left(P_1(R)+P_2(R)+2 R^2 P_2^2(R)-2 R^2 P_1(R) P_2(R)\right),\\
\label{E_2}
E_2=& 1 -12 P_1(R) + \frac{3}{2} P_2(R) + R^2 P_1(R)-\frac{1}{2} R^2 P_2(R)  +\frac{5}{2} R^2 P_2^2(R) -7 R^2 P_1(R) P_2(R) -\frac{1}{2} R^4 P_1(R) P_2^2(R), \\ 
\label{E_3}
E_3 =& - \frac{1}{2 R} \left( 24 P_1(R) - 4 R^2 P_1(R) +2 R^2 P_2(R)  -3 R^2 P_2^2(R) +22 R^2 P_1(R) P_2(R) + 2 R^4 P_1(R) P_2^2(R)\right).
\end{align}
\end{theorem}

In \cite{FriedmanHu2008}, Friedman and Hu transformed the stationary problem with Dirichlet boundary condition into the problem $\hat{F}(\widetilde{R}, \mu)=0$, where $\hat{F}(\ \cdot \  ,\mu)$ is a mapping between two Banach spaces. The computation of $\hat{\mu}_2'(0)$ depends on the second-order derivatives of $\hat{F}$, which are complicated due to the presence of the free boundary. The difficulty in determining the sign of $\hat{\mu}_2'(0)$ lies in the mixed calculations of Bessel functions $I_{1/2}$, $I_{3/2}$ and $I_{5/2}$. \eqref{cite dmu(0)} is a beautiful result obtained by the skillful inequality. In fact, the Dirichlet boundary condition corresponds to the case of infinite nutrient supply rate, i.e., $\beta= \infty$ in \eqref{sigma boundary}. For the case $0<\beta<\infty$, the introduction of $\beta$ adds difficulties to the process of calculating $\mu_2'(0)$ and its sign determination. Notice that $E_1$ in \eqref{dmu(0)} corresponds to the term in bracket of \eqref{cite dmu(0)}. $E_2$ and $E_3$ are new terms caused by the involvement of $\beta$. We need to study them further to determine the sign of $\mu_2'(0)$.

\begin{figure}[htbp]
    \centering
    \includegraphics[width=0.75\textwidth]{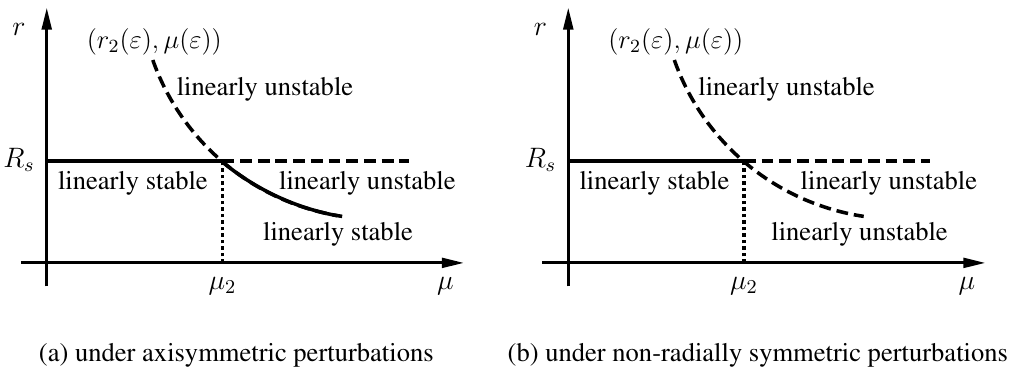}
    \caption{\centering Linear Stability}
    \label{figure3}
\end{figure}

The main result of this paper is the following theorem on the linear stability of the first bifurcation.
\begin{theorem}
\label{result2}
For problem \eqref{1.1}--\eqref{1.5} with $c=0$ and $\beta(t)  = \beta$, the axial symmetry stationary bifurcation solution $(\sigma_2(\varepsilon),p_2(\varepsilon),r_2(\varepsilon))$ is linearly unstable under non-radially symmetric perturbations when $|\varepsilon|$ is sufficiently small. (see Fig.\ref{figure3} (b))
\end{theorem}

\eqref{1.1}--\eqref{1.5} with $c=0$ and $\beta(t) = \beta$ is equivalent to $\frac{\partial \widetilde{R}}{\partial t} = -H(\widetilde{R},\mu)$(see \eqref{evolution equation}), where $H$ is given by \eqref{H}. Since the spectrum of $-H_{\widetilde{R}}(0, \mu_2)$ is in the left complex half-plane except the eigenvalue $0$(see Section \ref{sec4}), the linear stability of the stationary bifurcation solution $(\sigma_2(\varepsilon),p_2(\varepsilon), r_2(\varepsilon))$ is determined by the $0$-group eigenvalues of $-H_{\widetilde{R}}(\widetilde{R}_2(\varepsilon),\mu_2(\varepsilon))$($\widetilde{R}_2(\varepsilon)=r_2(\varepsilon)-R$), i.e., the totality of the
eigenvalues of $-H_{\widetilde{R}}(\widetilde{R}_2(\varepsilon),\mu_2(\varepsilon))$ arising from the splitting of the eigenvalue $0$
of the unperturbed operator $-H_{\widetilde{R}}(\widetilde{R}_2(0),\mu_2(0))$. In studying the linear stability of radially symmetric stationary solutions, the solutions of the linearized problem and the linearized operator $-H_{\widetilde{R}}(0,\mu)$ have the explicit expressions, which further implies that all eigenvalues of $-H_{\widetilde{R}}(0,\mu)$
admit the explicit forms as well. Since the stationary bifurcation solution $(\sigma_2(\varepsilon),p_2(\varepsilon),$ $ r_2(\varepsilon))$ is non-radially symmetric, 
we can't give the explicit expressions of the solutions for the linearized problem, while we only get the existence and uniqueness by partial differential equation theory. It follows that the expressions of $0$-group eigenvalues of $-H_{\widetilde{R}}(\widetilde{R}_2(\varepsilon),\mu_2(\varepsilon))$ can't be explicitly determined, making the analysis of the $0$-group difficult. By the principle of exchange of stability (cf. \cite[\S I.7]{Bifurcation2012}), Theorem \ref{result1} implies that the stationary bifurcation solution $(\sigma_2(\varepsilon),p_2(\varepsilon),r_2(\varepsilon))$ is linearly stable for $\varepsilon>0$ and linearly unstable for $\varepsilon<0$ under axisymmetric perturbations (see Fig.\ref{figure3} (a)). A natural problem is what the result will be for the linear stability of the stationary bifurcation solution $(\sigma_2(\varepsilon),p_2(\varepsilon),r_2(\varepsilon))$ under non-radially symmetric perturbations? We come up with a method to solve this problem with the help of the bifurcation curve’s structure (see Remark \ref{remark}). This method not only avoids a lot of calculations and the utilization of the explicit expressions of the solutions, but also generalizes the results on the linear stability of the stationary bifurcation solution $(\sigma_2(\varepsilon),p_2(\varepsilon),$ $ r_2(\varepsilon))$ from axisymmetric perturbations to non-radially symmetric perturbations.

The organization of this paper is as follows. In Section 2, we present some preliminaries which will be used in the following proof. We give the local bifurcation structure of the bifurcation curve $\{(r_2(\varepsilon), \mu_2(\varepsilon))\}$ in Section 3. In Section 4, we show the linear stability of the stationary bifurcation solution $(\sigma_2(\varepsilon),p_2(\varepsilon),r_2(\varepsilon))$ under non-radially symmetric perturbations.

\section{Preliminaries}\label{sec2}
In this section, we present some formulas in the 3-dimensional spherical coordinate system and recall some properties of the modified Bessel function $I_n$, the spherical harmonic functions $Y_{n,m}$ and $P_n$ defined in \eqref{P_n}. Moreover, we derive some new properties.

Let $\vec{e}_r, \vec{e}_\theta, \vec{e}_{\varphi}$ be the unit normal vectors in the $r, \theta, \varphi$ directions respectively, i.e.,
\begin{align*}
\vec{e}_r=(\sin \theta \cos \varphi, \sin \theta \sin \varphi, \cos \theta), \qquad \vec{e}_\theta=(\cos \theta \cos \varphi, \cos \theta \sin \varphi,-\sin \theta), \qquad \vec{e}_{\varphi}=(-\sin \varphi, \cos \varphi, 0).
\end{align*}
Then, the gradient is given by
\begin{align}
\label{gradient}
\nabla_x=\vec{e}_r \dfrac{\partial}{\partial r} + \vec{e}_\theta \frac{1}{r} \dfrac{\partial}{\partial \theta} + \vec{e}_\varphi \frac{1}{r \sin \theta} \dfrac{\partial}{\partial \varphi} = \vec{e}_r \dfrac{\partial}{\partial r} + \frac{1}{r} \nabla_\omega,
\end{align}
where $\nabla_\omega = \vec{e}_\theta  \frac{\partial}{\partial \theta} +\vec{e}_\varphi \frac{1}{ \sin \theta} \frac{\partial}{\partial \varphi} $, and the Laplace operator is
\begin{align*}
\Delta= \frac{\partial^2}{\partial r^2}+\frac{2}{r} \frac{\partial}{\partial r}+ \frac{1}{r^2 \sin \theta} \frac{\partial}{\partial \theta}\left(\sin \theta \frac{\partial}{\partial \theta}\right)+\frac{1}{ r^2 \sin ^2 \theta} \frac{\partial^2}{\partial \varphi^2} =\frac{\partial^2}{\partial r^2}+\frac{2}{r} \frac{\partial}{\partial r}+\frac{1}{r^2} \Delta_\omega,
\end{align*}
where $\Delta_\omega=\frac{1}{\sin \theta} \frac{\partial}{\partial \theta}\left(\sin \theta \frac{\partial}{\partial \theta}\right)+\frac{1}{\sin ^2 \theta} \frac{\partial^2}{\partial \varphi^2}$.

\cite[Lemma 4.1 and its proof]{FriedmanHu2008} showed that for the surface $r = R + \varepsilon \widetilde{R}(\theta,\varphi)$, the mean curvature is 
\begin{align}
\label{mean curvature}
\left.\kappa\right|_{r=R + \varepsilon \widetilde{R}}=\frac{1}{R}- \varepsilon \frac{1}{R^2}\left(\widetilde{R}+\frac{1}{2} \Delta_\omega \widetilde{R}\right)+ \varepsilon^2 \frac{1}{R^3}\left(\widetilde{R}^2+\widetilde{R} \Delta_\omega \widetilde{R} \right)+O\left(\varepsilon^3\right),
\end{align}
and the unit outer normal vector is expressed as
\begin{align*}
\left. \vec{n} \right|_{r=R + \varepsilon \widetilde{R}}=&\frac{1}{\sqrt{1+\left|\varepsilon \nabla_\omega \widetilde{R}\right|^2 \Big/(R+ \varepsilon \widetilde{R})^2}}\left(\vec{e}_r-\frac{\varepsilon \widetilde{R}_\theta}{R+ \varepsilon \widetilde{R}} \vec{e}_\theta-\frac{\varepsilon \widetilde{R}_\varphi}{(R+ \varepsilon \widetilde{R}) \sin \theta} \vec{e}_\varphi\right).
\end{align*}
By a further calculation, we have
\begin{align}
\label{normal vector}
\left. \vec{n} \right|_{r=R + \varepsilon \widetilde{R}}=\left[ 1- \varepsilon^2 \frac{1}{2 R^2}\left(\widetilde{R}_\theta^2+\frac{\widetilde{R}_\varphi^2}{ \sin ^2 \theta}\right)  \right] \vec{e}_r -\varepsilon \frac{ \widetilde{R}_\theta}{R} \vec{e}_\theta- \varepsilon \frac{ \widetilde{R}_{\varphi}}{R \sin \theta} \vec{e}_{\varphi}+O\left(\varepsilon^3\right).
\end{align}

Next, we introduce some properties about $I_n(r)$, where $n$ and $r$ are non-negative real numbers.

It is well-known that $I_n(r)$ satisfies
\begin{align}
\label{I'}
&I_{n }'(r) =  - \frac{n }{r} I_{n }(r) + I_{n-1}(r) ~ (n \geq 1),\qquad I_n'(r) = \frac{n}{r} I_n(r) + I_{n+1}(r).
\end{align}
By \eqref{I'}, we directly compute to get  
\begin{align}
\label{dr I_n+1/2(r) / r^1/2}
&\frac{d}{d r}\left( \frac{I_{n+1 / 2}(r)}{r^{1 / 2}} \right)=\frac{I_{n-1 / 2}(r)}{r^{1 / 2}} - \frac{n+1}{r} \frac{I_{n+1 / 2}(r)}{r^{ 1/2 }}  = 
 \frac{n}{r} \frac{I_{n+ 1/2}(r)}{r^{1 / 2}} +  \frac{I_{n+3 / 2}(r )}{r^{1 / 2}},\\
\label{dr2 I_n+1/2(r) / r^1/2}
&\frac{d^2}{d r^2}\left( \frac{I_{n+1 / 2}(r )}{r^{1 / 2}} \right) = - \frac{n}{r^2} \frac{I_{n+1 / 2}(r)}{r^{1 / 2}}  + \frac{n}{r} \left( \frac{n}{r} \frac{I_{n+1 / 2}(r)}{r^{1 / 2}} +\frac{I_{n+3 / 2}(r)}{r^{1 / 2}} \right)   +\frac{I_{n+1 / 2}(r)}{r^{1 / 2}} -\frac{n+2}{r} \frac{I_{n+3 / 2}(r)}{r^{1 / 2}}  \nonumber\\
&\qquad \qquad \qquad \quad \ \ = \left[ \frac{n(n-1)}{r^2}+1 \right] \frac{I_{n+1 / 2}(r)}{r^{1 / 2}}- \frac{2}{r} \frac{I_{n+3 / 2}(r)}{r^{1 / 2}} ,\\
\label{dr3 I_n+1/2(r) / r^1/2}
&\frac{d^3}{d r^3}\left( \frac{I_{n+1/2}(r)}{r^{1 / 2}} \right)= - \frac{2 n(n-1)}{r^3} \frac{I_{n+1 / 2}(r)}{r^{1 / 2}}  + \left[ \frac{n(n-1)}{r^2}+1 \right] \left( \frac{n}{r} \frac{I_{n+1 / 2}(r)}{r^{1 / 2}} + \frac{I_{n+3 / 2}(r)}{r^{1 / 2}}\right)   \nonumber \\
&\qquad \qquad \qquad \qquad \qquad +\frac{2}{r^2} \frac{I_{n+3 / 2}(r)}{r^{1 / 2}}  - \frac{2}{r} \left(\frac{I_{n+1 / 2}(r)}{r^{1 / 2}}-\frac{n+2}{r} \frac{I_{n+3 / 2}(r)}{r^{1 / 2}} \right)  \nonumber \\
&\qquad \qquad \qquad \quad \ \ =   \left\{ -\frac{2 n(n-1)}{r^3} + \left[\frac{n(n-1)}{r^2}+1\right] \frac{n}{r}  -\frac{2}{r}\right\} \frac{I_{n+1 / 2}(r)}{r^{1 / 2}}  \nonumber \\
&\qquad \qquad \qquad \qquad \qquad + \left[ \frac{n(n-1)}{r^2} + 1 +\frac{2}{r^2} + \frac{2(n+2)}{r^2} \right] \frac{I_{n+3 / 2}(r)}{r^{1 / 2}} \nonumber \\
&\qquad \qquad \qquad \quad \ \ =  \frac{n-2}{r} \left[ \frac{n(n-1)}{r^2}+1 \right] \frac{I_{n+1 / 2}(r)}{r^{1 / 2}} + \left( \frac{n^2+n+6}{r^2}+1 \right) \frac{I_{n+3 / 2}(r)}{r^{1 / 2}}.
\end{align}

Now, we recall some properties of $P_n(n=0,1,2,\cdots)$ given in \cite{FriedmanHu2008}
\begin{align} 
\label{Pn > Pn+1}
&P_n(r) > P_{n+1}(r) >0, \qquad r>0,  \\
\label{P_0}
&P_0(r)=\frac{r \cosh (r) -\sinh (r) }{r^2 \sinh (r)}, \\
\label{P_n P_n+1}
&r^2 P_n(r) P_{n+1}(r)+(2 n+3) P_n(r) =1.
\end{align}
From \eqref{P_0} and \eqref{P_n P_n+1}, we obtain
\begin{align}
\label{P_1}
&P_1(r)=\frac{r^2 \sinh (r) - 3 r \cosh (r) + 3 \sinh (r)}{r^2( r\cosh (r) -\sinh (r))},\\   
\label{P_2}
&P_2(r)=\frac{r^3\cosh (r) - 6 r^2 \sinh (r) + 15 r \cosh (r) - 15 \sinh (r)}{r^2\left(r^2 \sinh (r) - 3r \cosh (r) +3 \sinh (r)\right)}.
\end{align}

By \eqref{dr I_n+1/2(r) / r^1/2}--\eqref{dr3 I_n+1/2(r) / r^1/2} and \eqref{P_n}, we have
\begin{align}
\label{dr I_n+1/2(R) / R^1/2}
&\frac{R^{1 / 2}}{I_{n+1 / 2}(R)} \left. \frac{d}{d r}\left( \frac{I_{n+1 / 2}(r)}{r^{1 / 2}} \right) \right|_{r=R} 
=\frac{n}{R}+R P_n(R),\\
\label{dr2 I_n+1/2(R) / R^1/2}
&\frac{R^{1 / 2}}{I_{n+1 / 2}(R)} \left. \frac{d^2}{d r^2}\left( \frac{I_{n+1 / 2}(r)}{r^{1 / 2}} \right) \right|_{r=R}
=\frac{n(n-1)}{R^2}+1-2 P_n(R),\\
\label{dr3 I_n+1/2(R) / R^1/2}
&\frac{R^{1 / 2}}{I_{n+1 / 2}(R)} \left. \frac{d^3}{d r^3}\left( \frac{I_{n+1 / 2}(r)}{r^{1 / 2}} \right) \right|_{r=R} 
=\frac{n-2}{R}\left[\frac{n(n-1)}{R^2}+1\right]+ \left(\frac{n^2+n+6}{R^2}+1\right) R P_n(R).
\end{align}

The family of the spherical harmonic functions $\left\{Y_{n, m}\right\}$ is a complete orthogonal basis in $L^2(\mathbb{S}^2)$, which is given by
\begin{align*}
Y_{n,m}(\theta, \varphi)=\sqrt{\frac{(n-m)!}{(n+m)!} \frac{2 n+1}{4 \pi}} L_n^{m}(\cos \theta) \mathrm{e}^{\mathrm{i} m \varphi}, \quad m=-n, \cdots,  n,
\end{align*}
where 
\begin{align*}
L_n^m(x)=\frac{(-1)^m}{2^n n!} \left(1-x^2\right)^{\frac{m}{2}} \frac{\mathrm{d}^{n+m}}{\mathrm{~d} x^{n+m}}\left(x^2-1\right)^n, \qquad  m=-n, \cdots, n .
\end{align*}
$Y_{n,m}$ satisfies
\begin{align}
\label{laplace Ynm}
&-\Delta_\omega Y_{n, m}=n(n+1) Y_{n, m},\\
\label{solution1}
&-\Delta\left(r^n Y_{n, m}\right)=0, \\ 
\label{solution2}
&-\Delta\left(\frac{I_{n+1 / 2}(r)}{r^{1 / 2}} Y_{n, m}\right)+\left(\frac{I_{n+1 / 2}(r)}{r^{1 / 2}} Y_{n, m}\right)=0. 
\end{align}
By direct computation, we have
\begin{align}
\label{product1}
&\langle Y_{2,0} Y_{2,m}, Y_{1,l} \rangle = 0, \qquad \left \langle \frac{\partial Y_{2,0}}{\partial \theta} \frac{\partial Y_{2,m}}{\partial \theta}, Y_{1,l} \right\rangle = 0, \qquad -2 \leq m \leq 2, -1 \leq l \leq 1, \\
\label{product2}
&\langle Y_{2,0} Y_{2,m}, Y_{2,s} \rangle = 0, \qquad  \left \langle \frac{\partial Y_{2,0}}{\partial \theta} \frac{\partial Y_{2,m}}{\partial \theta}, Y_{2,s} \right\rangle = 0, \qquad m \neq s, \\
\label{product3}
&\langle Y_{2,0} Y_{2,m}, Y_{2,m} \rangle = \left\{\begin{array}{lll}
\dfrac{1}{7} \sqrt{\dfrac{5}{\pi}}, & m= 0, \\
\dfrac{1}{14} \sqrt{\dfrac{5}{\pi}}, & m= \pm 1, \\
-\dfrac{1}{7} \sqrt{\dfrac{5}{\pi}}, & m= \pm 2,
\end{array}\right., \qquad \left \langle \frac{\partial Y_{2,0}}{\partial \theta} \frac{\partial Y_{2,m}}{\partial \theta}, Y_{2,m} \right\rangle =  \left\{\begin{array}{lll}
\dfrac{3}{7} \sqrt{\dfrac{5}{\pi}}, & m= 0, \\
\dfrac{3}{14} \sqrt{\dfrac{5}{\pi}}, & m= \pm 1, \\
-\dfrac{3}{7} \sqrt{\dfrac{5}{\pi}}, & m= \pm 2.
\end{array}\right. 
\end{align}

\section{The computation of $\mu_2^{\prime}(0)$}\label{sec3}
In this section, we compute $\mu_2^{\prime}(0)$, in which $\mu_2(\varepsilon)$ comes from the bifurcation curve $\{ r_2(\varepsilon), \mu_2(\varepsilon) \}$ obtained by \cite[Theorem 1.1]{HuangZhangHu2017bifurcation} (or see Introduction).

First, we recall some results in \cite{FriedmanLam2015radial,HuangZhangHu2017bifurcation}. Friedman and Lam \cite[Theorem 3.1]{FriedmanLam2015radial} showed that for any $\beta>0$ and $0<\tilde{\sigma}<1$, there exists a unique radially symmetric solution $(\sigma_s,p_s,R)$ for problem \eqref{sigma equation}--\eqref{dn p=0}. \cite[(2.23), (2.25) and (2.26)]{HuangZhangHu2017bifurcation} gave the explicit expression of $(\sigma_s,p_s,R)$ as 
\begin{align}
\label{sigma_s(r)}
&\sigma_s(r)=\frac{\beta}{\beta+R P_0(R)}\frac{ R^{1/2} }{ I_{1 / 2 }\left(R\right)} \frac{ I_{1 / 2}(r)}{r^{ 1/2 } }, \\
\label{p_s(r)} 
&p_s(r)=-\mu \sigma_s(r) + \frac{1}{6} \mu \tilde{\sigma} (r^2 - R^2) + \frac{\gamma}{R}+\frac{\mu \beta}{\beta+R P_0(R)}, 
\end{align}
where $R$ is uniquely determined by the equation
\begin{align}
\label{R_s equation}
&\frac{ \beta  P_0(R)}{ \beta+R P_0(R) } =  \dfrac{\tilde{\sigma}}{3}. 
\end{align}

From \eqref{sigma_s(r)}, \eqref{dr I_n+1/2(R) / R^1/2}--\eqref{dr3 I_n+1/2(R) / R^1/2} and \eqref{P_n P_n+1} with $n=0$, we obtain
\begin{align}
&\left.\sigma_s \right|_{r=R} =\frac{\beta}{\beta+R P_0(R)}, \nonumber \\
\label{dr sigma_s(R)}
&\left.\frac{\partial \sigma_s}{\partial r}\right|_{r=R} =\frac{\beta P_0(R)}{\beta+RP_0(R)} R ,\\
\label{dr2 sigma_s(R)}
&\left.\frac{\partial^2 \sigma_s}{\partial r^2}\right|_{r=R} =\frac{\beta}{\beta+RP_0(R)}\left( 1-2 P_0(R) \right) =\frac{\beta P_0(R)}{\beta+RP_0(R)}  \left( 1 +  R^2P_1(R) \right),  \\
\label{dr3 sigma_s(R)}
&\left.\frac{\partial^3 \sigma_s}{\partial r^3}\right|_{r=R}  =\frac{\beta}{\beta+RP_0(R)} \left[-\frac{2}{R}+\left(\frac{6}{R^2}+1 \right) R P_0(R) \right]  =\frac{\beta  P_0(R)}{\beta+RP_0(R)} R \left(1-2 \frac{1-3 P_0(R)}{R^2 P_0(R)}\right) \nonumber \\
&\qquad \qquad \ \ =\frac{\beta P_0(R) }{\beta+RP_0(R)} R  \left(1-2 P_1(R)\right) .
\end{align}

By \eqref{p_s(r)}, \eqref{dr2 sigma_s(R)}--\eqref{dr3 sigma_s(R)} and \eqref{R_s equation}, we have
\begin{align}
\label{dr2 p_s(R)}
\left.\frac{\partial^2 p_s}{\partial r^2}\right|_{r=R}=& -\left.\mu \frac{\partial^2 \sigma_s}{\partial r^2}\right|_{r=R}+\frac{1}{3} \mu \tilde{\sigma}= -\frac{\mu \beta P_0(R)}{\beta+R P_0(R)}\left( 1 + R^2 P_1(R) \right)+\frac{\mu \beta P_0(R)}{\beta+R P_0(R)} \nonumber \\
=&-\frac{\mu \beta  P_0(R)}{\beta+R P_0(R)} R^2 P_1(R),\\
\label{dr3 p_s(R)}
\left.\frac{\partial^3 p_s}{\partial r^3}\right|_{r=R}=& -\mu \left.\frac{\partial^3 \sigma_s}{\partial r^3}\right|_{r=R} = \frac{ \mu \beta  P_0(R)}{\beta+RP_0(R)} R \left(2 P_1(R) -1\right).
\end{align}

Set
\begin{align*}
& X^{l+\alpha}=\left\{\widetilde{R} \in C^{l+\alpha}( \mathbb{S}^2 ): \widetilde{R} \text { is } \pi \text {-periodic in } \theta, 2 \pi \text {-periodic in } \varphi\right\}, \\
& X_2^{l+\alpha}= \text { closure of the linear space spanned by }\left\{Y_{j, 0}, j=0,2,4, \ldots\right\} \text { in } X^{l+\alpha} .
\end{align*}
Let $\Omega_{\widetilde{R}}$ be the domain with the boundary $\partial \Omega_{\widetilde{R}}=\{(r,\theta,\varphi):~r=R+\widetilde{R}(\theta,\varphi),~0\leq \theta<2\pi,~0\leq \varphi \leq \pi\}$.
For any $\widetilde{R} \in X^{4+\alpha}$, let $(\sigma, p)$ be the solution of the following problem
\begin{align}
\label{R sigma equation}
& -\Delta \sigma + \sigma = 0, \quad &&x \in \Omega_{\widetilde{R}}, \\
\label{R sigma boundary}
& \frac{\partial \sigma}{\partial \vec{n}}+\beta(\sigma-1)=0, \quad &&x \in \partial \Omega_{\widetilde{R}}, \\
\label{R p equation}
& -\Delta p=\mu(\sigma-\tilde{\sigma}), \quad &&x \in \Omega_{\widetilde{R}}, \\
\label{R p boundary}
& p= \gamma \kappa, \quad && x \in \partial \Omega_{\widetilde{R}}. 
\end{align}

Set $\Phi(x) = |x|-R-\widetilde{R}=r-R-\widetilde{R}$. Then $\left. \vec{n} \right|_{\partial \Omega_{\widetilde{R}}} = \left. \frac{\nabla \Phi}{\left|\nabla \Phi\right|}\right|_{\partial \Omega_{\widetilde{R}}}$. Define $H:X^{4+\alpha} \times \mathbb{R} \to X^{1+\alpha}$ by
\begin{align}
\label{H}
H(\widetilde{R}, \mu) = \left. \left( \nabla{p} \cdot \nabla \Phi \right) \right|_{\partial \Omega_{\widetilde{R}}}.
\end{align}
In \cite{HuangZhangHu2017bifurcation}, Huang, Zhang and Hu defined a mapping $F:X^{4+\alpha} \times \mathbb{R} \to X^{1+\alpha}$ as
\begin{align*}
F(\widetilde{R}, \mu)= \left.\frac{\partial p}{\partial n}\right|_{\partial \Omega_{\widetilde{R}}}. 
\end{align*}
Thus, $H(\widetilde{R}, \mu) =  \left| \nabla \Phi \right| \big|_{\partial \Omega_{\widetilde{R}}} F(\widetilde{R}, \mu)$. 
Similar to the statement in Cui and Zhuang \cite[\S 3]{cui2018bifurcation}, we get $H,F \in C^{\infty}\left(O_r \times \mathbb{R}, X^{1+\alpha} \right)$, where $O_r  \subset X^{4+\alpha}$ is a spherical neighborhood of $0$ with radius $r$. Then $(\sigma,p,R+\widetilde{R})$ is a solution of \eqref{sigma equation}--\eqref{dn p=0} if and only if $(\widetilde{R},\mu)$ is a solution of $H(\widetilde{R}, \mu) = 0$.

By \eqref{gradient}, we have $\left| \nabla \Phi \right| \big|_{\partial \Omega_{\widetilde{R}}} = \Big( 1+\frac{\left|\nabla_\omega  \widetilde{R}\right|^2}{(R+ \widetilde{R})^2} \Big)^{\frac{1}{2}}$. A direct computation yields
\begin{align*}
\Bigg( 1+\frac{\left|\nabla_\omega  (\varepsilon \widetilde{R})\right|^2}{(R+ \varepsilon \widetilde{R})^2} \Bigg)^{\frac{1}{2}} = 1 + O(\varepsilon^2).
\end{align*}
Together with the fact $F(0,\mu) =0$, we obtain
\begin{align*}
H( \varepsilon \widetilde{R}, \mu) =  \left| \nabla \Phi \right| \big|_{\partial \Omega_{\varepsilon \widetilde{R}}} F( \varepsilon \widetilde{R}, \mu) &= \left( 1 + O(\varepsilon^2) \right) \left( \varepsilon F_{\widetilde{R}}(0,\mu) \widetilde{R} + \frac{\varepsilon^2}{2} F_{\widetilde{R} \widetilde{R}}(0,\mu)[\widetilde{R},\widetilde{R}] + O(\varepsilon^3) \right) \\
&= \varepsilon F_{\widetilde{R}}(0,\mu) \widetilde{R} + \frac{\varepsilon^2}{2} F_{\widetilde{R} \widetilde{R}}(0,\mu)[\widetilde{R},\widetilde{R}] + O(\varepsilon^3).
\end{align*}
Hence
\begin{align}
\label{DRH=DRF}
H_{\widetilde{R}}(0,\mu) = F_{\widetilde{R}}(0,\mu), \qquad H_{\widetilde{R} \widetilde{R}}(0,\mu)[\widetilde{R},\widetilde{R}] = F_{\widetilde{R} \widetilde{R}}(0,\mu)[\widetilde{R},\widetilde{R}].   
\end{align}
Then \cite[(3.38), (3.45) and (3.46)]{HuangZhangHu2017bifurcation} implied
\begin{align}
\label{DR H}
H_{\widetilde{R}}(0, \mu) \left[ Y_{n, m} \right]=B_n \left( \mu_n-\mu \right) Y_{n, m}.
\end{align}
where
\begin{align}
&\mu_n = \gamma \frac{3 n(n-1)(n+2)}{2 \tilde{\sigma} R^4}  \frac{\beta+\dfrac{n}{R}+R P_n(R)}{(n+\beta R) P_1(R)-(1+\beta R) P_n(R)}, \nonumber \\
\label{B_n}
&B_n=\left\{\begin{array}{ll}
\gamma \dfrac{n(n-1)(n+2)}{2 R^3 } \dfrac{1}{\mu_n}, &\qquad  n>0,\\
-\dfrac{\tilde{\sigma} R}{3}  \dfrac{\beta R (P_0(R)-P_1(R)) + P_0(R)}{\beta+R P_0(R)}, &\qquad n=0.
\end{array}\right.
\end{align}
Moreover, \cite[Lemma 3.2]{HuangZhangHu2017bifurcation} showed
\begin{align}
\label{B_n>0}  
B_n\left\{\begin{array}{ll}
< 0, &\quad n=0, \\
=0, &\quad n=1, \\
>0, &\quad n \geqslant 2,
\end{array} \quad R \in(0,+\infty),\right.
\end{align}
and \cite[Lemma 3.3]{HuangZhangHu2017bifurcation} implied
\begin{align}
\label{increasing}
\mu_n \text{ is  strictly increasing for } n \geq 2.
\end{align}

Specially, for $n=2$ and using \eqref{R_s equation}, we get
\begin{align}
\mu_2 =& \gamma \frac{12}{ R^4} \frac{\beta+R P_0(R)}{3 \beta P_0(R)}  \frac{\beta+\dfrac{2}{R}+R P_2(R)}{R\left[\left(\beta+\dfrac{2}{R}+R P_2(R)\right) P_1(R)-\left(\beta+\dfrac{1}{R}+R P_1(R)\right) P_2(R)\right]} \nonumber\\
\label{mu_2}
=&\gamma \frac{4}{ R^5} \frac{\beta+R P_0(R)}{\beta  P_0(R)}  \frac{\beta+\dfrac{2}{R}+R P_2(R)}{\left(\beta+\dfrac{2}{R}+R P_2(R)\right) P_1(R)-\left(\beta+\dfrac{1}{R}+R P_1(R)\right) P_2(R)} .
\end{align}

Let $\widetilde{R}_n(\varepsilon) = r_n(\varepsilon) - R$. For $n=2$, $(\sigma_2(\varepsilon),p(\varepsilon),r_2(\varepsilon))$ is a bifurcation solution of \eqref{sigma equation}--\eqref{dn p=0} with $\mu=\mu_2(\varepsilon)$ (see Introduction), which leads to
\begin{align}
\label{bifurcation curve}
H(\widetilde{R}_2(\varepsilon),\mu_2(\varepsilon)) = 0,\qquad (\widetilde{R}_2(0),\mu_2(0)) = (0,\mu_2).
\end{align}
To determine the local bifurcation structure of the bifurcation curve $\{(r_2(\varepsilon),\mu_2(\varepsilon) )\}$, we need to compute $\mu_2^{\prime}(0)$. 
By \cite[(\uppercase\expandafter{\romannumeral1}.6.3)]{Bifurcation2012}, we get 
\begin{align}
\label{dmu equation}
\mu_2^{\prime}(0)  = -\frac{1}{2} \frac{\left\langle  H_{\widetilde{R} \widetilde{R}}\left(0, \mu_2\right) [Y_{2,0},Y_{2,0}], Y_{2,0}\right\rangle}{\langle H_{\widetilde{R} \mu}\left(0, \mu_2\right) [Y_{2,0}] , Y_{2,0}\rangle}.    
\end{align}
Next, we compute the second-order differential of $H$ involved in \eqref{dmu equation}.

Take $\widetilde{R} = \varepsilon Y_{2,0}$ in problem \eqref{R sigma equation}--\eqref{R p boundary} and consider the expansion of its solution $(\sigma,p)$ with the form
\begin{align}
\label{sigma expansion2} 
&\sigma=\sigma_s+\varepsilon \sigma_1+\varepsilon^2 \sigma_2+O\left(\varepsilon^3\right), \\
\label{p expansion2}  
&p=p_s+\varepsilon p_1+\varepsilon^2 p_2+O\left(\varepsilon^3\right).
\end{align}

\cite[(3.15)--(3.18)]{HuangZhangHu2017bifurcation} showed that $(\sigma_1,p_1)$ satisfies 
\begin{align}
\label{sigma_1 equation}
&-\Delta \sigma_1+\sigma_1=0, \quad && x \in  B_R, \\
\label{sigma_1 boundary}
&\frac{\partial \sigma_1}{\partial r}+\beta \sigma_1=-  \left.\left(\frac{\partial^2 \sigma_s}{\partial r^2}+\beta \frac{\partial \sigma_s}{\partial r}\right)\right|_{r=R} Y_{2,0}, \quad && x \in \partial B_R, \\
\label{p_1 equation}
&-\Delta p_1=\mu \sigma_1, \quad && x \in  B_R ,\\
\label{p_1 boundary}
&p_1=-\gamma\frac{1}{R^2}\left( Y_{2,0} + \frac{1}{2} \Delta_\omega Y_{2,0} \right), \quad && x \in \partial B_R,
\end{align}
and \cite[(3.24) and (3.26)]{HuangZhangHu2017bifurcation} implied 
\begin{align*}
\sigma_{1}(r, \theta)
=&\frac{- \beta \left( 1-2P_0(R)+\beta R P_0(R) \right) R^{ 1/2 }}{ \left( \beta+RP_0(R) \right) \left( \dfrac{2}{R} I_{5 / 2}(R)+I_{7 / 2}(R)+\beta I_{5 / 2}(R) \right)}  \frac{I_{5 / 2}(r)}{r^{ 1/2 }} Y_{2, 0}(\theta), \\
p_{1}(r, \theta)=&-\mu \sigma_1(r, \theta) +\left[ \gamma \frac{2}{R^2} -\frac{\mu \beta \left( 1-2P_0(R)+\beta R P_0(R) \right) }{ \left( \beta+RP_0(R) \right) \left( \beta+\dfrac{2}{R}+R P_2(R)  \right) } \right] \left(\frac{r}{R}\right)^2 Y_{2, 0}(\theta) .
\end{align*}
We use \eqref{P_n} and \eqref{P_n P_n+1} with $n=1$ to rewrite $(\sigma_{1},p_{1})$ as
\begin{align}
\label{sigma_1(r)}
\sigma_1(r, \theta)
=&- \frac{ \beta  P_0(R) }{\beta+RP_{0}(R) }
\frac{ \beta +  \dfrac{1}{R} +RP_1(R) }{ \beta+\dfrac{2}{R}+R P_2(R) } R \frac{R^{ 1/2 } }{ I_{5 / 2}\left(R\right)}
\frac{I_{5 / 2}(r)}{r^{ 1/2 }} Y_{2, 0}(\theta),\\
\label{p_1(r)}
p_1(r, \theta)=&-\mu \sigma_1(r, \theta) +\left( \gamma \frac{2}{R^2} -\frac{\mu \beta P_0(R)}{ \beta+RP_0(R) }\frac{\beta +  \dfrac{1}{R} +RP_1(R)}{ \beta+\dfrac{2}{R}+R P_2(R) } R \right) \left(\frac{r}{R}\right)^2 Y_{2, 0}(\theta) .
\end{align}

By \eqref{sigma_1(r)} and \eqref{dr I_n+1/2(R) / R^1/2}--\eqref{dr2 I_n+1/2(R) / R^1/2}, we get 
\begin{align}
\label{dr sigma_1(R)}
&\left.\frac{\partial \sigma_1}{\partial r}\right|_{r=R}=- \frac{ \beta  P_0(R) }{\beta+RP_{0}(R) }
\frac{\beta +  \dfrac{1}{R} +RP_1(R)}{ \beta+\dfrac{2}{R}+R P_2(R) } \left( 2 + R^2 P_2(R) \right) Y_{2,0},\\
\label{dtheta sigma_1(R)}
&\left.\frac{\partial \sigma_1}{\partial \theta}\right|_{r=R}=- \frac{ \beta  P_0(R) }{\beta+RP_{0}(R) }
\frac{ \beta +  \dfrac{1}{R} +RP_1(R) }{ \beta+\dfrac{2}{R}+R P_2(R) } R \frac{\partial Y_{2,0}}{\partial \theta} ,\\
\label{dr2 sigma_1(R)}
&\left.\frac{\partial^{2} \sigma_{1}}{\partial r^{2}}\right|_{r=R}=- \frac{ \beta  P_0(R) }{\beta+RP_{0}(R) }
\frac{\beta +  \dfrac{1}{R} +RP_1(R)}{ \beta+\dfrac{2}{R}+R P_2(R) } R \left( \frac{2}{R^2}+1-2P_2(R) \right)   Y_{2,0}.
\end{align}
Together with \eqref{p_1(r)}, we have
\begin{align}
\label{dr p_1(R)}
\left.\frac{\partial p_1}{\partial r}\right|_{r=R}=& \frac{ \mu \beta  P_0(R) }{\beta+RP_{0}(R) }
\frac{\beta +  \dfrac{1}{R} +RP_1(R)}{ \beta+\dfrac{2}{R}+R P_2(R) } \left( 2+R^2 P_2(R) \right) Y_{2,0} \nonumber\\
&\quad+ \left( \gamma \frac{2}{R^{2}} -\frac{ \mu \beta  P_0(R) }{\beta+RP_{0}(R) }
\frac{\beta +  \dfrac{1}{R} +RP_1(R)}{ \beta+\dfrac{2}{R}+R P_2(R) } R \right) \frac{2}{R} Y_{2,0} \nonumber \\
=&\left( \frac{\mu \beta  P_0(R)}{\beta+R P_0(R)} \frac{\beta +  \dfrac{1}{R} +RP_1(R)}{ \beta+\dfrac{2}{R}+R P_2(R) } R^2 P_2(R)+ \gamma \frac{4}{R^3} \right) Y_{2,0},  \\
\label{dtheta p_1(R)}
\left.\dfrac{\partial p_1}{\partial \theta}\right|_{r=R} =& \frac{ \mu \beta  P_0(R) }{\beta+RP_{0}(R) }
\frac{ \beta +  \dfrac{1}{R} +RP_1(R) }{ \beta+\dfrac{2}{R}+R P_2(R) } R \frac{\partial Y_{2,0}}{\partial \theta} + \left( \gamma \frac{2}{R^2}-\frac{\mu \beta  P_0(R)}{\beta+R P_0(R)} \frac{\beta+\dfrac{1}{R}+R P_1(R)}{\beta+\dfrac{2}{R}+R P_2(R)} R \right) \frac{\partial Y_{2,0}}{\partial \theta} \nonumber \\
=& \gamma \frac{2}{R^2} \frac{\partial Y_{2,0}}{\partial \theta},\\
\label{dr2 p_1(R)}
\left.\frac{\partial^{2} p_{1}}{\partial r^{2}}\right|_{r=R}
&=\frac{\mu \beta  P_0(R)}{\beta+R P_0(R)} \frac{\beta +  \dfrac{1}{R} +RP_1(R)}{ \beta+\dfrac{2}{R}+R P_2(R) } R \left(\frac{2}{R^2}+1-2 P_2(R)\right) Y_{2, 0} \nonumber \\
&\quad +\left( \gamma \frac{2}{R^2} -\frac{\mu \beta P_0(R)}{\beta+R P_0(R)} \frac{\beta +  \dfrac{1}{R} +RP_1(R)}{ \beta+\dfrac{2}{R}+R P_2(R) }  R \right) \frac{2}{R^2} Y_{2, 0} \nonumber \\
&=  \left[ \frac{\mu \beta  P_0(R)}{\beta+R P_0(R)} \frac{\beta +  \dfrac{1}{R} +RP_1(R)}{ \beta+\dfrac{2}{R}+R P_2(R) } R \left(1-2 P_2(R)\right)+ \gamma \frac{4}{R^4} \right] Y_{2, 0}.
\end{align}

\begin{lemma}
$(\sigma_2,p_2)$ satisfies
\begin{align}
\label{sigma_2 equation}
&-\Delta \sigma_2 + \sigma_2 =0, \qquad && x \in B_R,    \\
\label{sigma_2 boundary}
&\frac{\partial \sigma_2}{\partial r}+\beta \sigma_2 =-\frac{1}{2} \left. \left( \frac{\partial^3 \sigma_s}{\partial r^3} + \beta \frac{\partial^2 \sigma_s}{\partial r^2} \right) \right|_{r=R} Y_{2 ,0}^2  
+ \left.\frac{1}{2 R^2} \frac{\partial \sigma_s}{\partial r}\right|_{r=R}\left(\frac{\partial Y_{2,0}}{\partial \theta}\right)^2 \nonumber \\
&\qquad \qquad \qquad - \left.\left(\frac{\partial^2 \sigma_1}{\partial r^2}+\beta \frac{\partial \sigma_1}{\partial r}\right)\right|_{r=R} Y_{2 ,0} + \frac{1}{R^2} \left. \frac{\partial \sigma_1}{\partial \theta}\right|_{r=R} \frac{\partial Y_{2, 0}}{\partial \theta}, && x \in \partial B_R,\\
\label{p_2 equation}
&-\Delta p_2=\mu \sigma_2,  \qquad  && x \in  B_R,\\
\label{p_2 boundary}
&p_2 = \gamma \frac{1}{R^3}\left(Y_{2,0}^2+Y_{2,0} \Delta_\omega Y_{2,0} \right) -  \left.\frac{1}{2} \frac{\partial^2 p_s}{\partial r^2}\right|_{r=R} Y_{2,0}^2 -  \left.\frac{\partial p_1}{\partial r}\right|_{r=R} Y_{2,0}, && x \in \partial B_R.
\end{align}
\end{lemma}
\begin{proof}
By \eqref{R sigma equation}, \eqref{sigma expansion2}, \eqref{sigma_1 equation} and $-\Delta \sigma_s + \sigma_s =0$, we get \eqref{sigma_2 equation}.

From \eqref{R sigma boundary}, \eqref{sigma expansion2}, \eqref{gradient} and \eqref{normal vector}, we obtain
\begin{align*}
\beta  =&\left.\left(\frac{\partial \sigma}{\partial n}+\beta \sigma\right)\right|_{r=R+\varepsilon Y_{2, 0}} \\
=&\left.\nabla \sigma \cdot n\right|_{r=R+\varepsilon Y_{2, 0}}+\left.\beta \sigma\right|_{r=R+\varepsilon Y_{2, 0}} \\
=&\left.\left( \frac{\partial \sigma}{\partial r} \vec{e}_r + \frac{1}{R+\varepsilon Y_{2, 0}} \frac{\partial \sigma}{\partial \theta} \vec{e}_\theta + \frac{1}{\left(R+\varepsilon Y_{2, 0}\right) \sin \theta} \frac{\partial \sigma}{\partial \varphi} \vec{e}_{\varphi}\right)\right|_{r=R+\varepsilon Y_{2, 0}}\\
&\cdot \left\{\left[1-\varepsilon^2 \frac{1}{2 R^2}\left(\frac{\partial Y_{2, 0}}{\partial \theta}\right)^2\right] \vec{e}_r-\varepsilon \frac{1}{R} \frac{\partial Y_{2, 0}}{\partial \theta} \vec{e}_\theta+O\left(\varepsilon^3\right)\right\}+\left.\beta \sigma\right|_{r=R+\varepsilon Y_{2, 0}}\\
=&\left. \frac{\partial \sigma}{\partial r}\right|_{r=R+\varepsilon Y_{2, 0}} \left[1-\varepsilon^2 \frac{1}{2 R^2}\left(\frac{\partial Y_{2, 0}}{\partial \theta}\right)^2\right] - \left.\varepsilon  \frac{1}{R+\varepsilon Y_{2, 0}} \frac{\partial \sigma}{\partial \theta}\right|_{r=R+\varepsilon Y_{2, 0}}   \frac{1}{R} \frac{\partial Y_{2, 0}}{\partial \theta} + \left.\beta \sigma\right|_{r=R+\varepsilon Y_{2, 0}}+O\left(\varepsilon^3\right)\\
=&\left[\left.\frac{\partial \sigma_s}{\partial r}\right|_{r=R}+\left.\varepsilon \frac{\partial^2 \sigma_s}{\partial r^2}\right|_{r=R} Y_{2, 0}+\left.\varepsilon^2 \frac{1}{2} \frac{\partial^3 \sigma_s}{\partial r^3}\right|_{r=R} Y_{2, 0}^2 +\varepsilon\left(\left.\frac{\partial \sigma_1}{\partial r}\right|_{r=R}+\left.\varepsilon \frac{\partial^2 \sigma_1}{\partial r^2}\right|_{r=R} Y_{2,0}\right)+\left.\varepsilon^2 \frac{\partial \sigma_2}{\partial r}\right|_{r=R}\right]\\
& \qquad \cdot \left[1-\varepsilon^2 \frac{1}{2 R^2}\left(\frac{\partial Y_{2, 0}}{\partial \theta}\right)^2\right]  -\left.\varepsilon^2 \frac{1}{R^2}  \frac{\partial \sigma_1}{\partial \theta}\right|_{r=R} \frac{\partial Y_{2, 0}}{\partial \theta} \\
&+\beta\left[\left.\sigma_s\right|_{r=R}+\varepsilon \left. \frac{\partial \sigma_s}{\partial r} \right|_{r=R} Y_{2, 0}+\left.\varepsilon^2 \frac{1}{2} \frac{\partial^2 \sigma_s}{\partial r^2}\right|_{r=R} Y_{2, 0}^2 + \varepsilon \left(\left.\sigma_1\right|_{r=R}+\left.\varepsilon \frac{\partial \sigma_1}{\partial r}\right|_{r=R} Y_{2, 0}\right)+\left.\varepsilon^2 \sigma_2\right|_{r=R}\right]+O\left(\varepsilon^3\right)\\
=&\left.\left(\frac{\partial \sigma_s}{\partial r}+\beta \sigma_s\right)\right|_{r=R}+\varepsilon\left[\left.\left(\frac{\partial^2 \sigma_s}{\partial r^2}+\beta \frac{\partial \sigma_s}{\partial r}\right)\right|_{r=R} Y_{2, 0}+\left.\left(\frac{\partial \sigma_1}{\partial r}+\beta \sigma_1\right)\right|_{r=R}\right]\\
&+\varepsilon^2 \left[ \frac{1}{2} \left. \left( \frac{\partial^3 \sigma_s}{\partial r^3} + \beta \frac{\partial^2 \sigma_s}{\partial r^2} \right) \right|_{r=R} Y_{2 ,0}^2 -\left.\frac{1}{2 R^2} \frac{\partial \sigma_s}{\partial r}\right|_{r=R}\left(\frac{\partial Y_{2,0}}{\partial \theta}\right)^2 + \left. \left( \frac{\partial^2 \sigma_1}{\partial r^2} + \beta \frac{\partial \sigma_1}{\partial r} \right) \right|_{r=R} Y_{2, 0} \right.\\
&\left. -\left.\frac{1}{R^2} \frac{\partial \sigma_1}{\partial \theta}\right|_{r=R} \frac{\partial Y_{2, 0}}{\partial \theta} + \left. \left(\frac{\partial \sigma_2}{\partial r}+\beta \sigma_2\right) \right|_{r=R} \right] + O(\varepsilon^3).
\end{align*}
Using $\left. \left[ \frac{\partial \sigma_s}{\partial r}+\beta \left( \sigma_s - 1 \right) \right] \right|_{r=R}=0$ and \eqref{sigma_1 boundary}, we get the coefficient of $\varepsilon^0$ in the right side of above equality to be $\beta$ and $\varepsilon^1$ term to be $0$. It follows that \eqref{sigma_2 boundary} holds.

\eqref{R p equation}, \eqref{p expansion2}, \eqref{p_1 equation} and \eqref{sigma expansion2} imply \eqref{p_2 equation}.

By \eqref{p expansion2}, \eqref{p_1 boundary}, $\left.p_s\right|_{r=R}= \gamma \frac{1}{R}$ and $\left. \frac{\partial p_s}{\partial r}\right|_{r=R} = 0$, we obtain
\begin{align}
\label{p_2(R) compute1}
\left.p\right|_{r=R+\varepsilon Y_{2,0}} =&\left.p_s\right|_{r=R}+\left.\varepsilon \frac{\partial p_s}{\partial r}\right|_{r=R} Y_{2, 0}+\left.\varepsilon^2 \frac{1}{2} \frac{\partial^2 p_s}{\partial r^2}\right|_{r=R} Y_{2,0}^2 +\varepsilon\left(\left.p_1\right|_{r=R}+\left.\varepsilon \frac{\partial p_1}{\partial r}\right|_{r=R} Y_{2, 0}\right)+\left.\varepsilon^2 p_2\right|_{r=R}+O\left(\varepsilon^3\right) \nonumber \\
=& \gamma \frac{1}{R}- \varepsilon \gamma \frac{1}{R^2}\left(Y_{2,0}+\frac{1}{2} \Delta_\omega Y_{2,0}\right) + \varepsilon^2 \left( \left.\frac{1}{2} \frac{\partial^2 p_s}{\partial r^2}\right|_{r=R} Y_{2, 0}^2+\left.\frac{\partial p_1}{\partial r}\right|_{r=R} Y_{2, 0}+\left.p_2\right|_{r=R} \right) +O\left(\varepsilon^3\right).
\end{align}
On the other hand, \eqref{R p boundary} and \eqref{mean curvature} lead to
\begin{align}
\label{p_2(R) compute2}
\left.p\right|_{r=R+\varepsilon Y_{2,0}}= \gamma \left.\kappa\right|_{r=R + \varepsilon Y_{2,0}} =\gamma \frac{1}{R}- \varepsilon \gamma \frac{1}{R^2}\left(Y_{2,0}+\frac{1}{2} \Delta_\omega Y_{2,0}\right)+ \varepsilon^2 \gamma \frac{1}{R^3}\left(Y_{2,0}^2+Y_{2,0} \Delta_\omega Y_{2,0} \right)+O\left(\varepsilon^3\right) .
\end{align}
Then \eqref{p_2(R) compute1} and \eqref{p_2(R) compute2} yield \eqref{p_2 boundary}. We complete the proof. 
\end{proof}

\begin{lemma}
The following assertion holds
\begin{align}
\label{sigma_2 compute}
&\langle \sigma_2,Y_{2, 0} \rangle = D_1 \frac{R^{1 / 2}}{I_{5 / 2}(R)} \frac{I_{5 / 2}(r)}{r^{1 / 2}},\\    
\label{p_2 compute}
&\langle p_2,Y_{2, 0} \rangle = - \mu \langle \sigma_2,Y_{2, 0} \rangle + D_2 \left( \frac{r}{R} \right)^2 ,
\end{align}    
where 
\begin{align}
\label{D_1}
&D_1=\dfrac{1}{7} \sqrt{\dfrac{5}{\pi}}  \frac{ \beta P_0(R) }{\beta+RP_{0}(R) }
\frac{1}{ \left( \beta+\dfrac{2}{R}+R P_2(R) \right)^2} \nonumber \\
&\qquad \ \cdot \left\{ -\frac{1}{2} \left[R \left(1-2 P_1(R)\right) + \beta  \left(1+R^2 P_1(R)\right)   - \frac{3}{R} \right] \left( \beta+\dfrac{2}{R}+R P_2(R) \right) \right. \nonumber \\
&\qquad \qquad \left. +\left[ -\frac{1}{R} + R (1 - 2 P_2(R) )   +\beta \left( 2+R^2 P_2(R) \right)  \right] \left( \beta+\dfrac{1}{R}+R P_1(R) \right) \right\},\\
\label{D_2}
&D_2= \mu D_1 + \dfrac{1}{7} \sqrt{\dfrac{5}{\pi}} \left( -\gamma \frac{9}{R^3}+\frac{1}{2} \frac{\mu \beta  P_0(R)}{\beta+R P_0(R)} R^2 P_1(R)-\frac{\mu \beta  P_0(R)}{\beta+R P_0(R)} \frac{\beta+\dfrac{1}{R}+R P_1(R)}{\beta+\dfrac{2}{R}+R P_2(R)} R^2 P_2(R) \right).
\end{align}
\end{lemma}
\begin{proof}
By \eqref{sigma_2 equation} and \eqref{solution2}, $\langle \sigma_2,Y_{2, 0} \rangle$ has the form of \eqref{sigma_2 compute}. Together with \eqref{dr I_n+1/2(R) / R^1/2} and \eqref{P_n}, we get
\begin{align}
\label{D_1 compute}
\left\langle \left.\left(\frac{\partial \sigma_2}{\partial r}+\beta \sigma_2\right)\right|_{r=R},Y_{2, 0}  \right\rangle =D_1\left(\beta+\frac{2}{R}+R P_2(R)\right).
\end{align}
To determine $D_1$, we compute $\left\langle \left.\left(\frac{\partial \sigma_2}{\partial r}+\beta \sigma_2\right)\right|_{r=R},Y_{2, 0}  \right\rangle$ by \eqref{sigma_2 boundary}.

From \eqref{dr sigma_s(R)}--\eqref{dr3 sigma_s(R)} and  \eqref{dr sigma_1(R)}--\eqref{dr2 sigma_1(R)}, we have 
\begin{align*}
&-\left.\frac{1}{2}  \left(\frac{\partial^3 \sigma_s}{\partial r^3}+\beta \frac{\partial^2 \sigma_s}{\partial r^2}\right)\right|_{r=R} Y_{2,0}^2  =-\frac{1}{2} \frac{\beta P_0(R)}{\beta+R P_0(R)}\left[ R\left(1-2 P_1(R)\right) + \beta\left(1+R^2 P_1(R)\right)\right] Y_{2,0}^2,\\
&\frac{1}{2}  \frac{1}{R^2} \left. \frac{\partial \sigma_s}{\partial r}\right|_{r=R} \left( \frac{\partial Y_{2,0}}{\partial \theta} \right)^2 =\frac{1}{2}  \frac{\beta P_0(R)}{ \beta+R P_0(R) } \frac{1}{R} \left( \frac{\partial Y_{2,0}}{\partial \theta} \right)^2,\\
&-\left.\left(\frac{\partial^2 \sigma_1}{\partial r^2}+\beta \frac{\partial \sigma_1}{\partial r}\right)\right|_{r=R} Y_{2,0} = \frac{ \beta  P_0(R) }{\beta+RP_{0}(R) }
\frac{\beta +  \dfrac{1}{R} +RP_1(R)}{ \beta+\dfrac{2}{R}+R P_2(R) } \left[ R \left( \frac{2}{R^2}+1-2P_2(R) \right)   +\beta \left( 2+R^2 P_2(R) \right)  \right] Y_{2,0}^2,\\
&\left.\frac{1}{R^2}  \frac{\partial \sigma_1}{\partial \theta}\right|_{r=R} \frac{\partial Y_{2,0}}{\partial \theta}  = - \frac{ \beta P_0(R) }{\beta+RP_{0}(R) }
\frac{\beta +  \dfrac{1}{R} +RP_1(R)}{ \beta+\dfrac{2}{R}+R P_2(R) } \frac{1}{R} \left( \frac{\partial Y_{2,0}}{\partial \theta} \right)^2.
\end{align*}
Together with \eqref{sigma_2 boundary} and \eqref{product3}, we obtain 
\begin{align*}
&\left\langle \left.\left(\frac{\partial \sigma_2}{\partial r}+\beta \sigma_2\right)\right|_{r=R},Y_{2, 0}  \right\rangle \nonumber \\
&=-\frac{1}{2} \dfrac{1}{7} \sqrt{\dfrac{5}{\pi}} \frac{\beta P_0(R)}{\beta+R P_0(R)}\left[ R\left(1-2 P_1(R)\right) + \beta\left(1+R^2 P_1(R)\right) \right]   +\frac{1}{2} \dfrac{3}{7} \sqrt{\dfrac{5}{\pi}}   \frac{\beta P_0(R)}{ \beta+R P_0(R) } \frac{1}{R} \nonumber \\
&\quad + \dfrac{1}{7} \sqrt{\dfrac{5}{\pi}} \frac{ \beta  P_0(R) }{\beta+RP_{0}(R) }
\frac{\beta +  \dfrac{1}{R} +RP_1(R)}{ \beta+\dfrac{2}{R}+R P_2(R) } \left[ R \left( \frac{2}{R^2}+1-2P_2(R) \right)   +\beta \left( 2+R^2 P_2(R) \right)  \right] \nonumber \\
&\quad -\dfrac{3}{7} \sqrt{\dfrac{5}{\pi}}  \frac{ \beta P_0(R) }{\beta+RP_{0}(R) }
\frac{\beta +  \dfrac{1}{R} +RP_1(R)}{ \beta+\dfrac{2}{R}+R P_2(R) } \frac{1}{R}  \nonumber\\
&=-\frac{1}{2} \dfrac{1}{7} \sqrt{\dfrac{5}{\pi}}  \frac{\beta P_0(R)}{\beta+R P_0(R)}\left[ R \left(1-2 P_1(R)\right) + \beta \left(1+R^2 P_1(R)\right)  - \frac{3}{R} \right] \nonumber  \\
&\quad + \dfrac{1}{7} \sqrt{\dfrac{5}{\pi}}  \frac{ \beta P_0(R) }{\beta+RP_{0}(R) }
\frac{\beta +  \dfrac{1}{R} +RP_1(R)}{ \beta+\dfrac{2}{R}+R P_2(R) }
\left[ -\frac{1}{R} + R (1 - 2 P_2(R) )   +\beta  \left( 2+R^2 P_2(R) \right)  \right].
\end{align*}
Combining \eqref{D_1 compute}, we have  
\begin{align*}
D_1=&-\frac{1}{2} \dfrac{1}{7} \sqrt{\dfrac{5}{\pi}}  \frac{\beta P_0(R)}{\beta+R P_0(R)} \frac{1}{ \beta+\dfrac{2}{R}+R P_2(R) } \left[  R \left(1-2 P_1(R)\right) + \beta  \left(1+R^2 P_1(R)\right)  - \frac{3}{R} \right] \nonumber  \\
& + \dfrac{1}{7} \sqrt{\dfrac{5}{\pi}}  \frac{ \beta P_0(R) }{\beta+RP_{0}(R) }
\frac{\beta +  \dfrac{1}{R} +RP_1(R)}{ \left( \beta+\dfrac{2}{R}+R P_2(R) \right)^2}
\left[ -\frac{1}{R} + R(1 - 2 P_2(R) )   +\beta  \left( 2+R^2 P_2(R) \right)  \right]. 
\end{align*}
It follows that \eqref{D_1} holds.

By \eqref{sigma_2 equation} and \eqref{p_2 equation}, we get 
\begin{align*}
-\Delta \left( p_2+ \mu \sigma_2 \right)=0,  \qquad x \in  B_R.    
\end{align*}
Then \eqref{solution1} implies that $\langle p_2,Y_{2, 0} \rangle$ has the form of \eqref{p_2 compute}. Together with \eqref{sigma_2 compute}, we have
\begin{align}
\label{D_2 compute}
\left\langle\left. p_2\right|_{r=R}, Y_{2,0}\right\rangle =  -\mu D_1 + D_2 .
\end{align}

From \eqref{p_2 boundary}, \eqref{laplace Ynm}, \eqref{dr2 p_s(R)} and \eqref{dr p_1(R)}, we obtain
\begin{align*}
\left.p_2\right|_{r = R}= &\gamma \frac{1}{R^3}\left(Y_{2,0}^2-6 Y_{2,0}^2\right)+\frac{1}{2} \frac{\mu \beta  P_0(R)}{\beta+R P_0(R)} R^2 P_1(R) Y_{2,0}^2 -\left( \frac{\mu \beta P_0(R)}{\beta+R P_0(R)} \frac{\beta +  \dfrac{1}{R} +RP_1(R)}{ \beta+\dfrac{2}{R}+R P_2(R) } R^2 P_2(R)+ \gamma \frac{4}{R^3} \right) Y_{2,0}^2 \nonumber \\
= & \left( -\gamma \frac{9}{R^3}+\frac{1}{2} \frac{\mu \beta P_0(R)}{\beta+R P_0(R)} R^2 P_1(R)-\frac{\mu \beta  P_0(R)}{\beta+R P_0(R)} \frac{\beta+\dfrac{1}{R}+R P_1(R)}{\beta+\dfrac{2}{R}+R P_2(R)} R^2 P_2(R) \right) Y_{2,0}^2.
\end{align*}
Then using \eqref{product3}, we derive
\begin{align*}
&\left\langle\left. p_2\right|_{r=R}, Y_{2,0}\right\rangle  = \dfrac{1}{7} \sqrt{\dfrac{5}{\pi}} \left( -\gamma \frac{9}{R^3}+\frac{1}{2} \frac{\mu \beta  P_0(R)}{\beta+R P_0(R)} R^2 P_1(R)-\frac{\mu \beta  P_0(R)}{\beta+R P_0(R)} \frac{\beta+\dfrac{1}{R}+R P_1(R)}{\beta+\dfrac{2}{R}+R P_2(R)} R^2 P_2(R) \right).
\end{align*}
Together with \eqref{D_2 compute} and \eqref{D_1}, we get \eqref{D_2}. We complete the proof.  
\end{proof}

\eqref{sigma_2 compute}--\eqref{D_2} and \eqref{dr I_n+1/2(R) / R^1/2} lead to
\begin{align}
\left\langle \left.\frac{\partial \sigma_2}{\partial r}\right|_{r=R},Y_{2,0} \right\rangle  & =D_1\left(\frac{2}{R}+R P_2(R)\right),  \nonumber \\
\label{dr p_2(R)}
\left\langle \left.\frac{\partial p_2}{\partial r}\right|_{r=R},Y_{2,0} \right\rangle & =-\mu \left\langle \left.\frac{\partial \sigma_2}{\partial r}\right|_{r=R},Y_{2,0} \right\rangle + D_2 \frac{2}{R}   \nonumber \\
& =- \mu D_1 \left(\frac{2}{R}+R P_2(R)\right)+ \mu D_1 \frac{2}{R}  \nonumber \\
&\quad + \dfrac{1}{7} \sqrt{\dfrac{5}{\pi}} \left( - \gamma \frac{9}{R^3}+\frac{1}{2} \frac{\mu \beta  P_0(R)}{\beta+R P_0(R)} R^2 P_1(R)-\frac{\mu \beta  P_0(R)}{\beta+R P_0(R)} \frac{\beta+\dfrac{1}{R}+R P_1(R)}{\beta+\dfrac{2}{R}+R P_2(R)} R^2 P_2(R) \right) \frac{2}{R} \nonumber\\
&=\dfrac{1}{7} \sqrt{\dfrac{5}{\pi}}  \frac{ \mu \beta P_0(R) }{\beta+RP_{0}(R) }
\frac{1}{ \left( \beta+\dfrac{2}{R}+R P_2(R) \right)^2} R P_2(R) \nonumber \\
&\quad \cdot \left\{ \frac{1}{2} \left[ R \left(1-2 P_1(R)\right) + \beta  \left(1+R^2 P_1(R)\right)  - \frac{3}{R} \right] \left( \beta+\dfrac{2}{R}+R P_2(R) \right) \right. \nonumber \\
&\quad \qquad \left. -\left[ -\frac{1}{R} + R(1 - 2 P_2(R) )   +\beta  \left( 2+R^2 P_2(R) \right)  \right] \left( \beta+\dfrac{1}{R}+R P_1(R) \right) \right\} \nonumber \\
&\quad + \dfrac{1}{7} \sqrt{\dfrac{5}{\pi}} \left( - \gamma \frac{18}{R^4}+ \frac{\mu \beta P_0(R)}{\beta+R P_0(R)} R P_1(R)-2\frac{\mu \beta  P_0(R)}{\beta+R P_0(R)} \frac{\beta+\dfrac{1}{R}+R P_1(R)}{\beta+\dfrac{2}{R}+R P_2(R)} R P_2(R) \right).
\end{align}

Similar to the proof of \cite[Lemma 4.2]{FriedmanHu2008} and using \eqref{DRH=DRF}, we get
\begin{align}
\label{dmu compute1}
\frac{1}{2}  H_{\widetilde{R} \widetilde{R}}(0, \mu) [Y_{2,0},Y_{2,0}] &= \frac{1}{2}  F_{\widetilde{R} \widetilde{R}}(0, \mu) [Y_{2,0},Y_{2,0}] \nonumber \\
&=\left.\frac{1}{2} \frac{\partial^3 p_s}{\partial r^3}\right|_{r=R} Y_{2,0}^2 + \left. \frac{\partial^2 p_1}{\partial r^2}\right|_{r=R} Y_{2,0}  -\left.\frac{1}{R^2}  \frac{\partial p_1}{\partial \theta}\right|_{r=R} \frac{\partial Y_{2,0}}{\partial \theta}  + \left.\frac{\partial p_2}{\partial r}\right|_{r=R} .
\end{align}    
Hence \eqref{dmu compute1},  \eqref{dr3 p_s(R)}, \eqref{dr2 p_1(R)}, \eqref{dtheta p_1(R)},  \eqref{dr p_2(R)}  and \eqref{product3} imply
\begin{align}
\label{dmu compute2} 
&\left\langle \frac{1}{2} H_{\widetilde{R} \widetilde{R}}(0, \mu_2) [Y_{2,0},Y_{2,0}], Y_{2,0}\right\rangle \nonumber \\
&\quad= \frac{1}{2} \dfrac{1}{7} \sqrt{\dfrac{5}{\pi}} \frac{ \mu_2 \beta  P_0(R)}{\beta+RP_0(R)} R \left(2 P_1(R) -1\right) \nonumber \\
&\qquad  +\dfrac{1}{7} \sqrt{\dfrac{5}{\pi}} \left[ \frac{\mu_2 \beta  P_0(R)}{\beta+R P_0(R)} \frac{\beta +  \dfrac{1}{R} +RP_1(R)}{ \beta+\dfrac{2}{R}+R P_2(R) } R \left(1-2 P_2(R)\right)+ \gamma \frac{4}{R^4} \right] \nonumber \\
&\qquad  -\dfrac{3}{7} \sqrt{\dfrac{5}{\pi}} \gamma \frac{2}{R^4} \nonumber \\
&\qquad  +\dfrac{1}{7} \sqrt{\dfrac{5}{\pi}}  \frac{ \mu_2 \beta P_0(R) }{\beta+RP_{0}(R) }
\frac{1}{ \left( \beta+\dfrac{2}{R}+R P_2(R) \right)^2} R P_2(R) \nonumber \\
&\qquad \quad \cdot \left\{ \frac{1}{2} \left[ R \left(1-2 P_1(R)\right) + \beta  \left(1+R^2 P_1(R)\right)  - \frac{3}{R} \right] \left( \beta+\dfrac{2}{R}+R P_2(R) \right) \right. \nonumber \\
&\qquad \qquad \quad \left. -\left[ -\frac{1}{R} + R(1 - 2 P_2(R) )   +\beta  \left( 2+R^2 P_2(R) \right)  \right] \left( \beta+\dfrac{1}{R}+R P_1(R) \right) \right\}  \nonumber \\
&\qquad  + \dfrac{1}{7} \sqrt{\dfrac{5}{\pi}} \left( -\gamma \frac{18}{R^4}+ \frac{\mu_2 \beta P_0(R)}{\beta+R P_0(R)} R P_1(R)-2\frac{\mu_2 \beta  P_0(R)}{\beta+R P_0(R)} \frac{\beta+\dfrac{1}{R}+R P_1(R)}{\beta+\dfrac{2}{R}+R P_2(R)} R P_2(R) \right) \nonumber \\
&\quad=\dfrac{1}{7} \sqrt{\dfrac{5}{\pi}} \frac{ \mu_2 \beta P_0(R) }{\beta+RP_{0}(R) }
\frac{1}{ \left( \beta+\dfrac{2}{R}+R P_2(R) \right)^2} \nonumber \\
&\qquad \quad \cdot \left\{  \frac{1}{2} R \left(2 P_1(R) -1\right) \left( \beta+\dfrac{2}{R}+R P_2(R) \right)^2  \right. \nonumber \\
&\qquad \qquad \quad + R \left(1-2 P_2(R)\right) \left(\beta +  \dfrac{1}{R} +RP_1(R)\right) \left( \beta+\dfrac{2}{R}+R P_2(R) \right) \nonumber \\
&\qquad \qquad \quad  +\frac{1}{2} R P_2(R) \left[ R \left(1-2 P_1(R)\right) + \beta  \left(1+R^2 P_1(R)\right)  - \frac{3}{R} \right] \left( \beta+\dfrac{2}{R}+R P_2(R) \right) \nonumber \\
&\qquad \qquad \quad  -R P_2(R) \left[ -\frac{1}{R} + R(1 - 2 P_2(R) )   +\beta  \left( 2+R^2 P_2(R) \right)  \right] \left( \beta+\dfrac{1}{R}+R P_1(R) \right)   \nonumber \\
&\qquad \qquad \quad + R P_1(R) \left( \beta+\dfrac{2}{R}+R P_2(R) \right)^2 \nonumber \\
&\qquad \qquad \quad \left. - 2 R P_2(R) \left(\beta +  \dfrac{1}{R} +RP_1(R)\right) \left( \beta+\dfrac{2}{R}+R P_2(R) \right) \right\} \nonumber \\
&\qquad \quad - \dfrac{1}{7} \sqrt{\dfrac{5}{\pi}} \gamma \frac{20}{R^4}.
\end{align}

By \eqref{mu_2}, we obtain
\begin{align*}
-\gamma \frac{20}{R^4} =  -5R \frac{\mu_2 \beta  P_0(R) }{\beta+R P_0(R)}\cdot \frac{\left( \beta+ \dfrac{2}{R}+R P_2(R) \right) P_1(R)-\left( \beta+ \dfrac{1}{R}+R P_1(R) \right) P_2(R)}{\beta+ \dfrac{2}{R}+R P_2(R)}.    
\end{align*}
Together with \eqref{dmu compute2}, we derive
\begin{align}
\label{dmu compute3}
&\left\langle \frac{1}{2}  H_{\widetilde{R} \widetilde{R}}\left(0, \mu_2\right) [Y_{2,0},Y_{2,0}], Y_{2,0}\right\rangle \nonumber \\
&\quad =\frac{1}{7} \sqrt{\frac{{5}}{\pi}} \frac{ \mu_2 \beta P_0(R) }{\beta+RP_{0}(R) }
\frac{1}{ \left( \beta+\dfrac{2}{R}+R P_2(R) \right)^2} \nonumber \\
&\qquad \quad \cdot \left\{  \frac{1}{2} R \left(2 P_1(R) -1\right) \left( \beta+\dfrac{2}{R}+R P_2(R) \right)^2  \right. \nonumber \\
&\qquad \qquad \quad + R \left(1-2 P_2(R)\right) \left(\beta +  \dfrac{1}{R} +RP_1(R)\right) \left( \beta+\dfrac{2}{R}+R P_2(R) \right) \nonumber \\
&\qquad \qquad \quad  +\frac{1}{2} R P_2(R) \left[ R \left(1-2 P_1(R)\right) + \beta  \left(1+R^2 P_1(R)\right)  - \frac{3}{R} \right] \left( \beta+\dfrac{2}{R}+R P_2(R) \right) \nonumber \\
&\qquad \qquad \quad  -R P_2(R) \left[ -\frac{1}{R} + R(1 - 2 P_2(R) )   +\beta  \left( 2+R^2 P_2(R) \right)  \right] \left( \beta+\dfrac{1}{R}+R P_1(R) \right)   \nonumber \\
&\qquad \qquad \quad + R P_1(R) \left( \beta+\dfrac{2}{R}+R P_2(R) \right)^2 \nonumber \\
&\qquad \qquad \quad  - 2 R P_2(R) \left(\beta +  \dfrac{1}{R} +RP_1(R)\right) \left( \beta+\dfrac{2}{R}+R P_2(R) \right)  \nonumber \\
&\qquad \qquad \quad  -5 R P_1(R)  \left( \beta+\dfrac{2}{R}+R P_2(R) \right)^2  \nonumber \\
&\qquad \qquad \quad \left. +5 R P_2(R) \left(\beta +  \dfrac{1}{R} +RP_1(R)\right) \left( \beta+\dfrac{2}{R}+R P_2(R) \right) \right\} \nonumber \\
&\quad =\frac{1}{7} \sqrt{\frac{{5}}{\pi}} \frac{ \mu_2 \beta P_0(R) }{\beta+RP_{0}(R) }
\frac{1}{ \left( \beta+\dfrac{2}{R}+R P_2(R) \right)^2} \nonumber \\
&\qquad \quad \cdot \left\{ -R \left( \frac{1}{2} + 3  P_1(R)\right)\left(\beta+\frac{2}{R}+R P_2(R)\right)^2 \right.\nonumber \\
&\qquad \qquad \quad + R \left( 1 + P_2(R) \right)\left(\beta+\frac{1}{R}+R P_1(R)\right)\left(\beta+\frac{2}{R}+R P_2(R)\right)  \nonumber \\
& \qquad \qquad \quad + \frac{1}{2} R P_2(R) \left[ R \left(1-2 P_1(R)\right) + \beta  \left(1+R^2 P_1(R)\right)  - \frac{3}{R} \right] \left( \beta+\dfrac{2}{R}+R P_2(R) \right)  \nonumber \\
&\qquad \qquad \quad \left. -R P_2(R) \left[ -\frac{1}{R} + R(1 - 2 P_2(R) )   +\beta  \left( 2+R^2 P_2(R) \right)  \right] \left( \beta+\dfrac{1}{R}+R P_1(R) \right) \right\}  \nonumber \\
&\quad =\frac{1}{7} \sqrt{\frac{{5}}{\pi}} \frac{ \mu_2 \beta P_0(R) }{\beta+RP_{0}(R) }
\frac{1}{ \left( \beta+\dfrac{2}{R}+R P_2(R) \right)^2} \left( E_1 R \beta^2 + E_2 \beta + E_3 \right),
\end{align}
where 
\begin{align}
\label{E_1 compute}
E_1= -\frac{1}{2} \left(-1 + 6P_1(R) + P_2(R) + 2R^2 P_2^2(R) - R^2 P_1(R) P_2(R)    \right),    
\end{align}
and $E_2$ and $E_3$ are given by \eqref{E_2} and \eqref{E_3}, respectively.
Using \eqref{P_n P_n+1} with $n=1$, \eqref{E_1 compute} is rewritten as \eqref{E_1}.

By \cite[Lemma 4.6]{FriedmanHu2008}, we get $E_1<0$. We next show $E_2<0$ and $E_3<0$.

From \eqref{E_2} and \eqref{P_n P_n+1} with $n=1$, we obtain
\begin{align*}
E_2=& -7 P_1(R) + \frac{3}{2} P_2(R) + R^2 P_1(R)-\frac{1}{2} R^2 P_2(R)  +\frac{5}{2} R^2 P_2^2(R) -6 R^2 P_1(R) P_2(R) -\frac{1}{2} R^4 P_1(R) P_2^2(R) \nonumber \\ 
=& \left( -7 + R^2 \right) P_1(R) + \left( \frac{3}{2} -\frac{1}{2} R^2 \right) P_2(R)  +\frac{5}{2} R^2 P_2^2(R) -6 R^2 P_1(R) P_2(R) -\frac{1}{2} R^4 P_1(R) P_2^2(R). 
\end{align*}
Together with \eqref{P_1} and \eqref{P_2}, we have
\begin{align}
\label{E_2 compute}
E_2=& \left( -7 + R^2 \right) \frac{ R^2 \sinh (R) - 3 R \cosh (R) + 3 \sinh (R)}{R^2( R\cosh (R) - \sinh (R))} \nonumber\\ 
&+ \left( \frac{3}{2} -\frac{1}{2} R^2 \right) \frac{R^3\cosh (R) - 6 R^2 \sinh (R) + 15 R \cosh (R) - 15 \sinh (R)}{R^2\left(R^2 \sinh (R) - 3R \cosh (R) +3 \sinh (R)\right)} \nonumber \\
&+ \frac{5}{2} R^2 \left( \frac{R^3\cosh (R) - 6 R^2 \sinh (R) + 15 R \cosh (R) - 15 \sinh (R)}{R^2\left(R^2 \sinh (R) - 3R \cosh (R) +3 \sinh (R)\right)} \right)^2 \nonumber \\
&-6 R^2 \frac{ R^2 \sinh (R) - 3 R \cosh (R) + 3 \sinh (R)}{R^2( R\cosh (R) - \sinh (R))} \frac{R^3\cosh (R) - 6 R^2 \sinh (R) + 15 R \cosh (R) - 15 \sinh (R)}{R^2\left(R^2 \sinh (R) - 3R \cosh (R) +3 \sinh (R)\right)} \nonumber \\
&-\frac{1}{2} R^4 \frac{ R^2 \sinh (R) - 3 R \cosh (R) + 3 \sinh (R)}{R^2( R\cosh (R) - \sinh (R))} \left( \frac{R^3\cosh (R) - 6 R^2 \sinh (R) + 15 R \cosh (R) - 15 \sinh (R)}{R^2\left(R^2 \sinh (R) - 3R \cosh (R) +3 \sinh (R)\right)} \right)^2 \nonumber \\
=&\frac{G_1(R)}{8 R^2 \left(-3 R \cosh(R) + (3 + R^2) \sinh(R) \right)^2 (R \cosh(R) - \sinh(R))},
\end{align}
where 
\begin{align}
\label{G_1}
G_1(R)=& (-1269 + 741 R^2 + 366 R^4 + 44 R^6) R \cosh(R) \nonumber \\
&+ 3 (423 + 317 R^2 + 30 R^4) R \cosh(3R) \nonumber \\ 
&+  \left( 1269 - 648 R^2 - 606 R^4 - 146 R^6 - 8 R^8 \right)  \sinh(R) \nonumber\\
&-  \left( 423 + 1476 R^2 + 378 R^4 + 10 R^6 \right) \sinh(3 R).  
\end{align}
The result of the simplification of $E_2$ has been verified with matlab.

Notice the denominator of \eqref{E_2 compute} is greater than $0$ for $R>0$. To prove $E_2<0$, it is sufficient to prove $G_1(R)<0$. Through the following well-known formulas
\begin{align*}
\sinh(x)=\sum_{n=0}^{\infty} \frac{x^{2 n+1}}{(2 n+1)!},  \qquad \cosh(x)=\sum_{n=0}^{\infty} \frac{x^{2 n}}{(2 n)!}, \qquad x \in(-\infty,+\infty),   
\end{align*}
we rewrite \eqref{G_1} as
\begin{align*}
G_1(R)= & \left( -1269 + 741 R^2 + 366 R^4 + 44 R^6 \right) R \sum_{n=0}^{\infty} \frac{1}{(2 n)!} R^{2 n} \nonumber \\
&+ 3 \left( 423 + 317 R^2 + 30 R^4 \right) R \sum_{n=0}^{\infty} \frac{3^{2 n}}{(2 n)!} R^{2 n}\nonumber \\ 
& +\left( 1269 - 648 R^2 - 606 R^4 - 146 R^6 - 8 R^8 \right)  \sum_{n=0}^{\infty} \frac{1}{(2 n+1)!} R^{2 n+1} \nonumber\\
&-  \left( 423 + 1476 R^2 + 378 R^4 + 10 R^6 \right) \sum_{n=0}^{\infty} \frac{3^{2 n+1}}{(2 n+1)!} R^{2 n+1}  \\
=& - \sum_{n=0}^{\infty} \frac{1269}{(2 n)!} R^{2 n + 1} + \sum_{n=0}^{\infty} \frac{741}{(2 n)!} R^{2 n + 3}  + \sum_{n=0}^{\infty} \frac{366}{(2 n)!} R^{2 n + 5} + \sum_{n=0}^{\infty} \frac{44}{(2 n)!} R^{2 n + 7} \nonumber \\
&+ \sum_{n=0}^{\infty} \frac{423 \cdot 3^{2n+1}}{(2 n)!} R^{2 n + 1} + \sum_{n=0}^{\infty} \frac{317 \cdot 3^{2n+1} }{(2 n)!} R^{2 n + 3}  + \sum_{n=0}^{\infty} \frac{30 \cdot 3^{2n+1}}{(2 n)!} R^{2 n + 5}  \nonumber \\
&+\sum_{n=0}^{\infty} \frac{1269}{(2 n +1)!} R^{2 n + 1} - \sum_{n=0}^{\infty} \frac{648}{(2 n +1)!} R^{2 n + 3}  - \sum_{n=0}^{\infty} \frac{606}{(2 n +1)!} R^{2 n + 5} - \sum_{n=0}^{\infty} \frac{146}{(2 n +1)!} R^{2 n + 7} \nonumber \\
&- \sum_{n=0}^{\infty} \frac{8}{(2 n +1)!} R^{2 n + 9} \nonumber \\
&-\sum_{n=0}^{\infty} \frac{423 \cdot 3^{2n+1}}{(2 n+1)!} R^{2 n + 1} - \sum_{n=0}^{\infty} \frac{1476 \cdot 3^{2n+1} }{(2 n +1)!} R^{2 n + 3}  - \sum_{n=0}^{\infty} \frac{378 \cdot 3^{2n+1}}{(2 n+1)!} R^{2 n + 5} - \sum_{n=0}^{\infty} \frac{10 \cdot 3^{2n+1}}{(2 n+1)!} R^{2 n + 7} \nonumber \\
=& - \sum_{n=0}^{\infty} \frac{1269}{(2 n)!} R^{2 n + 1} + \sum_{n=1}^{\infty} \frac{741}{(2 n-2)!} R^{2 n + 1}  + \sum_{n=2}^{\infty} \frac{366}{(2 n-4)!} R^{2 n + 1} + \sum_{n=3}^{\infty} \frac{44}{(2 n-6)!} R^{2 n + 1} \nonumber \\
&+ \sum_{n=0}^{\infty} \frac{423 \cdot 3^{2n+1}}{(2 n)!} R^{2 n + 1} + \sum_{n=1}^{\infty} \frac{317 \cdot 3^{2n-1} }{(2 n-2)!} R^{2 n + 1}  + \sum_{n=2}^{\infty} \frac{30 \cdot 3^{2n-3}}{(2 n-4)!} R^{2 n + 1}  \nonumber \\
&+\sum_{n=0}^{\infty} \frac{1269}{(2 n +1)!} R^{2 n + 1} - \sum_{n=1}^{\infty} \frac{648}{(2 n -1)!} R^{2 n + 1}  - \sum_{n=2}^{\infty} \frac{606}{(2 n -3)!} R^{2 n + 1} - \sum_{n=3}^{\infty} \frac{146}{(2 n -5)!} R^{2 n + 1} \nonumber \\
&- \sum_{n=4}^{\infty} \frac{8}{(2 n -7)!} R^{2 n + 1} \nonumber \\
&-\sum_{n=0}^{\infty} \frac{423 \cdot 3^{2n+1}}{(2 n+1)!} R^{2 n + 1} - \sum_{n=1}^{\infty} \frac{1476 \cdot 3^{2n-1} }{(2 n -1)!} R^{2 n + 1}  - \sum_{n=2}^{\infty} \frac{378 \cdot 3^{2n-3}}{(2 n-3)!} R^{2 n + 1} - \sum_{n=3}^{\infty} \frac{10 \cdot 3^{2n-5}}{(2 n-5)!} R^{2 n + 1} \nonumber \\
=&-1269 R - \frac{1269}{2!} R^3 - \frac{1269}{4!} R^5  - \frac{1269}{6!} R^7  + 741 R^3 + \frac{741}{2!} R^5 + \frac{741}{4!} R^7  + 366 R^5 + \frac{366}{2!} R^7 + 44 R^7  \nonumber \\
&+ 423 \cdot 3 R + \frac{423 \cdot 3^3}{2!} R^3 + \frac{423 \cdot 3^5}{4!} R^5 + \frac{423 \cdot 3^7}{6!} R^7  + 317 \cdot 3 R^3 + \frac{317 \cdot 3^3}{2!} R^5 + \frac{317 \cdot 3^5}{4!} R^7   \nonumber\\
&+ 30 \cdot 3 R^5 + \frac{30 \cdot 3^3}{2!} R^7 + 1269 R + \frac{1269}{3!} R^3 + \frac{1269}{5!} R^5  + \frac{1269}{7!} R^7  - 648 R^3 - \frac{648}{3!} R^5 -  \frac{648}{5!} R^7  \nonumber \\
&- 606 R^5 - \frac{606}{3!} R^7  - 146 R^7 - 423 \cdot 3 R - \frac{423 \cdot 3^3}{3!} R^3 - \frac{423 \cdot 3^5}{5!} R^5 -\frac{423 \cdot 3^7}{7!} R^7  \nonumber \\
&- 1476 \cdot 3 R^3 - \frac{1476 \cdot 3^3}{3!} R^5 - \frac{1476 \cdot 3^5}{5!} R^7  - 378 \cdot 3 R^5 - \frac{378 \cdot 3^3}{3!} R^7  - 10 \cdot 3 R^7  \nonumber \\
&+ \sum_{n=4}^{\infty} \frac{R^{2 n + 1}}{(2 n +1)!} \nonumber \\
&\quad \cdot \left[ -1269(2 n +1) + 741 (2 n +1) (2 n) (2 n -1) + 366 (2 n +1) (2 n) (2 n -1) (2 n -2) (2 n -3) \right. \nonumber \\
&\qquad
+ 44 (2 n +1) (2 n) (2 n -1) (2 n -2) (2 n -3) (2 n -4) (2 n -5) \nonumber \\
&\qquad 
+ 423 (2 n +1) 3^{2 n+1} + 317 (2 n +1) (2 n) (2 n -1) 3^{2 n-1} \nonumber \\
&\qquad + 30 (2 n +1) (2 n) (2 n -1) (2 n -2) (2 n -3) 3^{2 n-3} + 1269 -648 (2 n +1) (2 n) \nonumber \\
&\qquad  - 606 (2 n +1) (2 n) (2 n -1) (2 n -2)  - 146  (2 n +1) (2 n) (2 n -1) (2 n -2) (2 n -3) (2 n -4)  \nonumber \\
&\qquad  - 8 (2 n +1) (2 n) (2 n -1) (2 n -2) (2 n -3) (2 n -4) (2 n -5) (2 n -6)   \nonumber \\
&\qquad  -423 \cdot 3^{2n+1} - 1476 (2 n +1) (2 n) 3^{2n-1} - 378 (2 n +1) (2 n) (2n-1) (2n-2) 3^{2n-3} \nonumber \\
&\qquad  \left. -10 (2 n +1) (2 n) (2n-1) (2n-2) (2n-3) (2n-4) 3^{2n-5} \right] \nonumber \\
=& \left( -1269 + 423 \cdot 3 + 1269 - 423 \cdot 3 \right) R  \nonumber \\
&+ \left(  -\frac{1269}{2!} + 741 + \frac{423 \cdot 3^3}{2!} + 317 \cdot 3 + \frac{1269}{3!} - 648 - \frac{423 \cdot 3^3}{3!} - 1476 \cdot 3 \right) R^3 \nonumber \\
&+ \left( - \frac{1269}{4!} + \frac{741}{2!} + 366 + \frac{423 \cdot 3^5}{4!} + \frac{317 \cdot 3^3}{2!} + 30 \cdot 3 + \frac{1269}{5!} - \frac{648}{3!} \right. \nonumber \\
&\qquad \left. - 606 - \frac{423 \cdot 3^5}{5!} - \frac{1476 \cdot 3^3}{3!} - 378 \cdot 3 \right) R^5 \nonumber \\
&+ \left( - \frac{1269}{6!} + \frac{741}{4!} + \frac{366}{2!} + 44 + \frac{423 \cdot 3^7}{6!} + \frac{317 \cdot 3^5}{4!} + \frac{30 \cdot 3^3}{2!} + \frac{1269}{7!} -  \frac{648}{5!} - \frac{606}{3!}  \right.  \nonumber \\
& \qquad \left. - 146 -\frac{423 \cdot 3^7}{7!} - \frac{1476 \cdot 3^5}{5!} - \frac{378 \cdot 3^3}{3!} - 10 \cdot 3 \right) R^7 \nonumber \\
&+ \sum_{n=4}^{\infty} \frac{R^{2 n + 1}}{(2 n +1)!} \nonumber \\
&\quad \cdot \left\{ -1269(2 n +1) + 741 (2 n +1) (2 n) (2 n -1) + 366 (2 n +1) (2 n) (2 n -1) (2 n -2) (2 n -3) \right. \nonumber \\
&\qquad
+ 44 (2 n +1) (2 n) (2 n -1) (2 n -2) (2 n -3) (2 n -4) (2 n -5)  + 1269 -648 (2 n +1) (2 n)  \nonumber \\
&\qquad  - 606 (2 n +1) (2 n) (2 n -1) (2 n -2)  - 146  (2 n +1) (2 n) (2 n -1) (2 n -2) (2 n -3) (2 n -4)  \nonumber \\
&\qquad - 8 (2 n +1) (2 n) (2 n -1) (2 n -2) (2 n -3) (2 n -4) (2 n -5) (2 n -6)   \nonumber \\
&\qquad 
+  \left[ 423 \cdot 3^6(2 n +1)  + 317 \cdot 3^4 (2 n +1) (2 n) (2 n -1)  \right. \nonumber \\
&\qquad \qquad + 30 \cdot 3^2 (2 n +1) (2 n) (2 n -1) (2 n -2) (2 n -3)   \nonumber \\
&\qquad \qquad -423 \cdot 3^6 - 1476 \cdot 3^4 (2 n +1) (2 n) - 378 \cdot 3^2 (2 n +1) (2 n) (2n-1) (2n-2) \nonumber \\
&\qquad \qquad \left. \left. -10 (2 n +1) (2 n) (2n-1) (2n-2) (2n-3) (2n-4) \right] 3^{2 n-5} \right\} \nonumber \\
=&-4 \sum_{n=4}^{\infty} \frac{n R^{2 n + 1}}{(2 n +1)!} \left( a_n + b_n  3^{2n-5} \right),
\end{align*}
where
\begin{align*}
&a_n=10305 - 20712n - 18498n^{2} + 77432n^{3} - 73920n^{4} + 31904n^{5} - 6528n^{6} + 512n^{7},\\
&b_n=-77235 + 114544n - 67962n^{2} + 20008n^{3} - 2880n^{4} + 160n^{5}.
\end{align*}
Then
\begin{align*}
&a_n=10305 + 20712n \left( -1 +n^2  \right) + \left(- 18498 + 56720n \right) n^{2} + \left( - 73920 + 31904n \right) n^{4} + \left( - 6528 + 512n \right) n^{6},\\
&b_n=\left( -77235 + 114544n \right) + \left(- 67962 + 20008n \right) n^{2} + \left( - 2880 + 160n \right) n^{4}.
\end{align*}
Notice that each bracket term on the right side of the equalities above is greater than $0$ when $n \geq 18$. By direct calculation, we obtain that $a_n + b_n  3^{2n-5} \geq 0$ for $4 \leq n \leq 17$ (see Table \ref{Tab}). It follows that $G_1(R)<0$ for $R>0$. Moreover, we get $E_2<0$.

\begin{table}[h]
\centering
\caption{Computation of $a_n + b_n 3^{2n-5}$ for $4\leq n \leq 17$}\label{Tab}
\begin{tabular}{*{11}{c}} 
\toprule
$n$ & $a_n + b_n 3^{2n-5}$ &  & $n$  & $a_n + b_n 3^{2n-5}$ &  & $n$  & $a_n + b_n 3^{2n-5}$ &  & $n$  & $a_n + b_n 3^{2n-5}$ \\ 
\midrule
4 & 0                    &  & 8  & 32037857280          &  & 12 & 7180490438922240     &  & 16 & 345732732222060165120  \\
5 & 0                    &  & 9  & 935872045056         &  & 13 & 115430374226534400   &  & 17 & 4609555109351370768384 \\
6 & 0                    &  & 10 & 21237572689920       &  & 14 & 1743676765861109760  &  &    &                        \\
7 & 656308224            &  & 11 & 412227655004160      &  & 15 & 25059810408839424000 &  &    &                        \\
\bottomrule
\end{tabular}
\end{table}

It remains to prove that $E_3 < 0$. Multiplying \eqref{P_n P_n+1} with $n=2$ by $-4R^2P_1(R)$ and the one with $n=1$ by $2R^2P_2(R)$, we get 
\begin{align*}
&-4 R^2 P_1(R)=-4 R^4 P_1(R) P_2(R) P_3(R)-28 R^2 P_1(R) P_2(R),\\
&2 R^2 P_2(R)=2 R^4 P_1(R) P_2^2(R)+10 R^2 P_1(R) P_2(R)
.
\end{align*}
Substituting the above equations into \eqref{E_3}, we have
\begin{align*}
E_3 =& - \frac{1}{2R} \left( 24 P_1(R) -4 R^4 P_1(R) P_2(R) P_3(R)-28 R^2 P_1(R) P_2(R)  + 2 R^4 P_1(R) P_2^2(R) \right. \nonumber \\
&\qquad \qquad \left. +10 R^2 P_1(R) P_2(R) -3 R^2 P_2^2(R) +22 R^2 P_1(R) P_2(R) + 2 R^4 P_1(R) P_2^2(R)\right) \nonumber \\
=& - \frac{1}{2R} \left[ 24 P_1(R) + 4 R^4 P_1(R) P_2(R) ( P_2(R) - P_3(R) )  + R^2 P_1(R) P_2(R) +3 R^2 P_2(R) ( P_1(R) - P_2(R) )  \right]. 
\end{align*}
\eqref{Pn > Pn+1} leads to $P_1(R) - P_2(R)>0$ and $P_2(R) - P_3(R)>0$ for $R>0$. Therefore, $E_3<0$.
So far, we have proved $E_1<0$, $E_2<0$ and $E_3<0$.

By \eqref{DR H}, we have $H_{\widetilde{R} \mu}(0, \mu_2) [Y_{2, 0}]= -B_2 Y_{2, 0}$. Together with \eqref{B_n}, we obtain 
\begin{align}
\label{dmu compute4}
\left\langle H_{\widetilde{R} \mu}(0, \mu_2) [Y_{2, 0}], Y_{2,0}\right\rangle =-B_2=-\gamma\frac{4}{R^3} \frac{1}{\mu_2}.
\end{align}

From \eqref{dmu equation}, \eqref{dmu compute3}, \eqref{dmu compute4}, $E_1<0$, $E_2<0$ and $E_3<0$, we derive \eqref{dmu(0)}.

\section{Linear stability}\label{sec4}
In this section, we study the linear stability of the stationary bifurcation solution $(\sigma_2(\varepsilon),p_2(\varepsilon),r_2(\varepsilon))$ under non-radially symmetric perturbations.

For any $\widetilde{R}(\ \cdot \ ,t) \in X^{4+\alpha}$, consider the domain $\Omega_{\widetilde{R}(t)}$ with the boundary $\partial \Omega_{\widetilde{R}(t)}: r=R+\widetilde{R}(\theta,\varphi,t)$.
Let $(\sigma,p)$ be the solution of problem \eqref{1.1}--\eqref{1.4} with $c=0$, $\beta(t) = \beta$ and $\Omega(t) = \Omega_{\widetilde{R}(t)}$. Since 
\begin{align*}
V_n = \Bigg(1+\frac{\left|\nabla_\omega \widetilde{R}\right|^2}{(R+\widetilde{R})^2}\Bigg)^{-1 / 2} \frac{\partial \widetilde{R}}{\partial t},    
\end{align*}
problem \eqref{1.1}--\eqref{1.5} with $c=0$ and $\beta(t) = \beta$ is equivalent to 
\begin{align}
\label{evolution equation}
\dfrac{\partial \widetilde{R}}{\partial t} = -H(\widetilde{R},\mu),
\end{align}
where $H$ is given by \eqref{H}.

Set $R_*$ as an equilibrium of \eqref{evolution equation} for the parameter $\mu_*$, i.e., $H(R_*,\mu_*)=0$. Notice problem \eqref{1.1}--\eqref{1.5} with $c=0$ and $\beta(t) = \beta$ is invariant under the translation. Similar to the statement in Pan and Xing \cite[\S 3.4]{PanXingBifurcation2022}, we get
\begin{align}
\label{DRH Y1m = 0}
H_{\widetilde{R}}(R_*,\mu_*) Y_{1,m} = 0, \qquad  -1 \leq m \leq 1.   
\end{align}
The translation of the known solution isn't a new one, which leads us to define the generalized principle of linearized stability as follow  
\begin{align}
\label{principle}
\begin{array}{ll}
&\text{the equilibrum } R_* \text{ is linearly stable for the parameter $\mu_*$ if }\\
&\operatorname{Ker}\left[-H_{\widetilde{R}}(R_*,\mu_* )\right] = \operatorname{span}\left\{ Y_{1,0}, \dfrac{1}{2}(Y_{1,-1}-Y_{1,1}), \dfrac{i}{2}(Y_{1,-1}+Y_{1,1}) \right\}\\
&\text{and the spectrum of } -H_{\widetilde{R}}(R_*,\mu_* ) \text{ is in the left complex half-plane except $0$.} 
\end{array}
\end{align}

By \eqref{bifurcation curve}, $\widetilde{R}_2(\varepsilon)$ is an equilibrium of \eqref{evolution equation} for the parameter $\mu_2(\varepsilon)$. We next study the linear stability of $(\widetilde{R}_2(\varepsilon) ,\mu_2(\varepsilon))$.

From \eqref{DR H}, we obtain that all spectral points of $-H_{\widetilde{R}}(0, \mu_2)$ are real. We now show that if $\lambda \neq -B_n \left(\mu_n-\mu_2 \right)$, $\lambda I+H_{\widetilde{R}}(0, \mu_2): X^{4+\alpha} \to X^{1+\alpha}$ is a bijection, i.e, the spectrum of $-H_{\widetilde{R}}(0, \mu_2)$ consists entirely of eigenvalues. In fact, by \eqref{DR H} and the continuity of $-H_{\widetilde{R}}(0, \mu_2)$, we obtain that if $\lambda \neq -B_n \left(\mu_n-\mu_2 \right)$, $\lambda I+H_{\widetilde{R}}(0, \mu_2)$ is an injection. Similar to the proof of \cite[Theorem 7.3]{WellEscher2011}, we get that $-H_{\widetilde{R}}(0, \mu_2)$ is a Fredholm operator of index $0$. Since $I:X^{4+\alpha} \to X^{1+\alpha}$ is compact, $\lambda I+H_{\widetilde{R}}(0, \mu_2)$ is also a Fredholm operator of index $0$. Together with the fact that $\lambda I+H_{\widetilde{R}}(0, \mu_2)$ is an injection, we derive that $\lambda I+H_{\widetilde{R}}(0, \mu_2)$ is a bijection. 

Combining \eqref{B_n>0}--\eqref{increasing}, we derive that the spectrum of $-H_{\widetilde{R}}(0, \mu_2)$ are in the left complex half-plane except $0$. 
Thus, to study the linear stability of $(\widetilde{R}_2(\varepsilon) ,\mu_2(\varepsilon))$, it is sufficient to study the $0$-group eigenvalues of $-H_{\widetilde{R}}(\widetilde{R}_2(\varepsilon),\mu_2(\varepsilon))$.

At first, we calculate the total number of the $0$-group eigenvalues of $H_{\widetilde{R}}(\widetilde{R}_2(\varepsilon),\mu_2(\varepsilon))$ when $|\varepsilon|$ is small.

\begin{lemma}
\label{8 eigenvalues}
There exists a positive real number $\delta$ such that for any $\varepsilon \in (-\delta,\delta)$, the number of
all eigenvalues in the $0$-group of $H_{\widetilde{R}}(\widetilde{R}_2(\varepsilon),\mu_2(\varepsilon))$ (counting multiplicity)
is $8$. 
\end{lemma}
\begin{proof}
Notice $H_{\widetilde{R}}(\widetilde{R}_2(0),\mu_2(0))=H_{\widetilde{R}}(0,\mu_2 )$. \eqref{DR H}, \eqref{B_n>0} and \eqref{increasing} imply that $0$ is an isolated eigenvalue of $H_{\widetilde{R}}(0,\mu_2 )$ with the eigenspace
\begin{align*}
\operatorname{Ker}\left[H_{\widetilde{R}}(0,\mu_2 )\right] = \operatorname{span}\left\{ Y_{1,0}, Y_{1,1}, Y_{1,-1}, Y_{2,0}, Y_{2,1}, Y_{2,-1},
Y_{2,2}, Y_{2,-2} \right\},
\end{align*}
i.e., $\operatorname{dim} \operatorname{Ker}\left[H_{\widetilde{R}}(0,\mu_2 )\right] = 8$.
Moreover, $X^{1+\alpha}$ have the following decomposition
\begin{align}
\label{decomposition2}
X^{1+\alpha} = \operatorname{Ker}\left[H_{\widetilde{R}}(0,\mu_2 )\right] \bigoplus  \operatorname{Im}\left[H_{\widetilde{R}}(0,\mu_2 )\right] .
\end{align}

In order to show that the algebraic multiplicity of the eigenvalue $0$ of $H_{\widetilde{R}}(0,\mu_2 )$ is equal to the geometric multiplicity $8$, it is sufficient to show that $\operatorname{Ker}\left[H_{\widetilde{R}}(0,\mu_2 )\right]^2 = \operatorname{Ker}\left[H_{\widetilde{R}}(0,\mu_2 )\right]$. 
Note $\operatorname{Ker}\left[H_{\widetilde{R}}(0,\mu_2 )\right]^2 \supset \operatorname{Ker}\left[H_{\widetilde{R}}(0,\mu_2 )\right]$. We next prove $\operatorname{Ker}\left[H_{\widetilde{R}}(0,\mu_2 )\right]^2 \subset \operatorname{Ker}\left[H_{\widetilde{R}}(0,\mu_2 )\right]$.

Assume $\xi  \in \operatorname{Ker}\left[H_{\widetilde{R}}(0,\mu_2 )\right]^2$. Then $H_{\widetilde{R}}(0,\mu_2 ) [\xi] \in \operatorname{Ker}\left[H_{\widetilde{R}}(0,\mu_2 )\right] \cap  \operatorname{Im}\left[H_{\widetilde{R}}(0,\mu_2 )\right]$. By \eqref{decomposition2}, we get $H_{\widetilde{R}}(0,\mu_2 ) [\xi] =0$, i.e., $\xi \in \operatorname{Ker}\left[H_{\widetilde{R}}(0,\mu_2 )\right]$.

\cite[II.5.1 and III.6.4]{Kato1976Perturbation} implied that for any small $|\varepsilon|$, the number of
all eigenvalues in the $0$-group of $H_{\widetilde{R}}(\widetilde{R}_2(\varepsilon),\mu_2(\varepsilon))$ (counting multiplicity) is equal to the algebraic multiplicity $8$ of the eigenvalue $0$ of $H_{\widetilde{R}}(0,\mu_2 )$. We get the result.
\end{proof}

Next, we use the following Proposition to prove that $(\widetilde{R}_2(\varepsilon),\mu_2(\varepsilon))$ is linearly unstable.

\begin{proposition}
\label{d eigenvalue}
There exist two eigenvalues $\lambda_1(\varepsilon)$ and $\lambda_2(\varepsilon)$ in the $0$-group of $H_{\widetilde{R}}(\widetilde{R}_2(\varepsilon),\mu_2(\varepsilon))$ such that $\lambda_1'(0) = \frac{1}{2}a$ and $\lambda_2'(0) = -\frac{3}{2}a$, where 
\begin{align}
\label{a}
a = -2 \mu_2'(0) \langle H_{\widetilde{R} \mu}(0,\mu_2) [Y_{2, 0}],Y_{2, 0} \rangle = \langle H_{\widetilde{R} \widetilde{R}}(0,\mu_2) \left[ Y_{2,0}, Y_{2, 0} \right],Y_{2, 0} \rangle  < 0 \ (\text{see \eqref{dmu equation} and \eqref{dmu compute3}}).  
\end{align}
\end{proposition}

\begin{remark}
\label{remark}
In fact, $a < 0$ plays an essential role in the proof of Theorem \ref{result2}, which is derived through the analysis of the bifurcation curve’s structure.    
\end{remark}

Before proving Proposition \ref{d eigenvalue}, we need the following lemma.

\begin{lemma}
\label{8 eigenvalues 1}
Define $Q_0$ by
\begin{align}
\label{Q_0}
Q_0 = -\frac{1}{2 \pi i} \int_{\Gamma} \left[ H_{\widetilde{R}}(0,\mu_2) - \zeta I \right]^{-1} d \zeta,
\end{align}
where $\Gamma$ is a simple closed curve such that all the eigenvalues are outside $\Gamma$ except $0$.
Then $Q_0 \{ H_{\widetilde{R} \widetilde{R}}(0, $ $\mu_2)  \left[ Y_{2,0}, \ \cdot \ \right] + H_{\widetilde{R} \mu}(0,\mu_2) \mu_2'(0) \} Q_0: \operatorname{Ker}\left[H_{\widetilde{R}}(0,\mu_2 )\right] \rightarrow \operatorname{Ker}\left[H_{\widetilde{R}}(0,\mu_2 )\right]$ satisfies
\begin{align*}
&Q_0 \{ H_{\widetilde{R} \widetilde{R}}(0,\mu_2) \left[ Y_{2,0}, \ \cdot \ \right] + H_{\widetilde{R} \mu}(0,\mu_2) \mu_2'(0) \} Q_0 \left[\begin{array}{l}
Y_{1,0} \ \ Y_{1,1} \ \ Y_{1,-1} \ \ Y_{2,0} \ \ Y_{2,1} \ \ Y_{2,-1} \ \
Y_{2,2} \ \ Y_{2,-2}
\end{array}\right]  \nonumber\\
&\qquad \quad =\left[\begin{array}{l}
Y_{1,0} \ \ Y_{1,1} \ \ Y_{1,-1} \ \ Y_{2,0} \ \ Y_{2,1} \ \ Y_{2,-1} \ \
Y_{2,2} \ \ Y_{2,-2}
\end{array}\right]  \left[\begin{array}{llllllll}
0 & & & & & & & \\
& 0 & & & & & & \\
& & 0 & & & & & \\
& & & \dfrac{1}{2}a & & & & \\
& & & & 0 & & \\
& & & & & 0 & \\
& & & & & & -\dfrac{3}{2}a\\
& & & & & & & -\dfrac{3}{2}a\\
\end{array}\right],
\end{align*}
where $a$ is given by \eqref{a}.
\end{lemma}
\begin{proof}
Since $0$ is an isolated eigenvalue of $H_{\widetilde{R}}(0,\mu_2)$, $Q_0$ is well-defined.

By \cite[III-\S 6.4, Theorem 6.17 and its proof]{Kato1976Perturbation}, $Q_0$ is a eigenprojection. We now show that $Q_0 X^{1+\alpha} = \operatorname{Ker}\left[H_{\widetilde{R}}(0,\mu_2 )\right]$. Let $\xi \in \operatorname{Ker}\left[H_{\widetilde{R}}(0,\mu_2 )\right]$. Then $H_{\widetilde{R}}(0,\mu_2 ) [\xi] =0$, so that
\begin{align*}
Q_0 \xi &=-\frac{1}{2 \pi i} \int_{\Gamma} \left[ H_{\widetilde{R}}(0,\mu_2) - \zeta I \right]^{-1} \xi \ d \zeta  = -\frac{1}{2 \pi i} \int_{\Gamma} \left[ H_{\widetilde{R}}(0,\mu_2) - \zeta I \right]^{-1} \left[ H_{\widetilde{R}}(0,\mu_2) - \zeta I \right] \left( -\dfrac{1}{\zeta} \xi \right) \ d \zeta \\
&= \frac{1}{2 \pi i} \int_{\Gamma} \dfrac{1}{\zeta} \xi \ d \zeta = \xi.   
\end{align*}
Therefore, $\operatorname{Ker}\left[H_{\widetilde{R}}(0,\mu_2 )\right] \subset Q_0 X^{1+\alpha}$. By \cite[III-\S6.5]{Kato1976Perturbation}, we obtain that the algebraic multiplicity of the eigenvalue $0$ of $H_{\widetilde{R}}(0,\mu_2 )$ is equal to $\operatorname{dim} Q_0 X^{1+\alpha}$. In the proof of Lemma \ref{8 eigenvalues}, we show that the algebraic multiplicity of the eigenvalue $0$ of $H_{\widetilde{R}}(0,\mu_2 )$ is equal to the geometric multiplicity $\operatorname{dim} \operatorname{Ker}\left[H_{\widetilde{R}}(0,\mu_2 )\right]=8$. Then $\operatorname{dim} Q_0 X^{1+\alpha} = \operatorname{dim} \operatorname{Ker}\left[H_{\widetilde{R}}(0,\mu_2 )\right]=8$. Hence $Q_0 X^{1+\alpha} = \operatorname{Ker}\left[H_{\widetilde{R}}(0,\mu_2 )\right]$.

\eqref{DRH Y1m = 0} implies $H_{\widetilde{R}}(\widetilde{R}_2(\varepsilon),\mu_2(\varepsilon)) [Y_{1,l}] = 0$ for $-1 \leq l \leq 1$. Taking the derivative at $\varepsilon=0$ yields
\begin{align*}
H_{\widetilde{R} \widetilde{R}}(0,\mu_2) \left[ Y_{2,0}, Y_{1, l} \right] + H_{\widetilde{R} \mu}(0,\mu_2) \mu_2'(0) [Y_{1, l}] = 0,
\end{align*}
which implies 
\begin{align}
\label{under Y1m}
Q_0 \{ H_{\widetilde{R} \widetilde{R}}(0,\mu_2) \left[ Y_{2,0}, \ \cdot \ \right] + H_{\widetilde{R} \mu}(0,\mu_2) \mu_2'(0) \} Q_0 [Y_{1, l}] = 0, \qquad -1 \leq l \leq 1.  
\end{align}

To compute $H_{\widetilde{R} \widetilde{R}}(0,\mu_2) \left[ Y_{2,0}, Y_{2, m} \right]$, we collect the $\varepsilon \eta$-term in the expansion of $H(\varepsilon Y_{2,0} + \eta Y_{2, m},\mu_2) (m\not=0)$. Taking $\widetilde{R} = \varepsilon Y_{2,0} + \eta Y_{2, m}$, we consider the expansion of its solution $(\sigma,p)$ with the form
\begin{align*}
&\sigma = \sigma_s + \varepsilon \sigma_1 + \eta w_1 + \varepsilon^2 \sigma_2 + \varepsilon \eta w_2 + \eta^2 w_3 + O\left( (\varepsilon^2+\eta^2)^{\frac{3}{2}} \right), \\  
&p = p_s + \varepsilon p_1 + \eta q_1 + \varepsilon^2 p_2 + \varepsilon \eta q_2 + \eta^2 q_3 + O\left( (\varepsilon^2+\eta^2)^{\frac{3}{2}} \right).
\end{align*}
Then 
\begin{align}
\label{DRR H Y20 Y2m}
H_{\widetilde{R} \widetilde{R}}\left(0, \mu_2\right)\left[Y_{2,0}, Y_{2, m}\right] =&\left.\frac{\partial^3 p_s}{\partial r^3}\right|_{r=R} Y_{2,0} Y_{2, m}+ \left. \frac{\partial^2 p_1}{\partial r^2}\right|_{r = R} Y_{2, m} -\left. \frac{1}{R^2} \frac{\partial p_1}{\partial \theta}\right|_{r= R} \frac{\partial Y_{2, m}}{\partial \theta} \nonumber \\
&+ \left. \frac{\partial^2 q_1}{\partial r^2}\right|_{r = R} Y_{2, 0} -\left. \frac{1}{R^2} \frac{\partial q_1}{\partial \theta}\right|_{r= R} \frac{\partial Y_{2, 0}}{\partial \theta} + \left. \frac{\partial q_2}{\partial r}\right|_{r = R},
\end{align}
where $p_s$ is given by \eqref{p_s(r)}, $p_1$ by \eqref{p_1(r)} and $q_1$ and $q_2$ satisfy the follow problem 
\begin{align}
\label{w_1 equation}
&-\Delta w_1+w_1=0, \quad && x \in  B_R, \\
\label{w_1 boundary}
&\frac{\partial w_1}{\partial r}+\beta w_1=-  \left.\left(\frac{\partial^2 \sigma_s}{\partial r^2}+\beta \frac{\partial \sigma_s}{\partial r}\right)\right|_{r=R} Y_{2,m}, \quad && x \in \partial B_R, \\
\label{q_1 equation}
&-\Delta q_1=\mu w_1, \quad && x \in  B_R ,\\
\label{q_1 boundary}
&q_1=-\gamma \frac{1}{R^2}\left( Y_{2,m} + \frac{1}{2} \Delta_\omega Y_{2,m} \right), \quad && x \in \partial B_R, \\
\label{w_2 equation}
&-\Delta w_2+w_2=0, \quad && x \in  B_R, \\
\label{w_2 boundary}
&\frac{\partial w_2}{\partial r}+\beta w_2= - \left.\left(\frac{\partial^3 \sigma_s}{\partial r^3}+\beta \frac{\partial^2 \sigma_s}{\partial r^2}\right)\right|_{r=R} Y_{2,0} Y_{2,m} + \left. \frac{1}{R^2} \frac{\partial \sigma_s}{\partial r}\right|_{r= R} \frac{\partial Y_{2, 0}}{\partial \theta} \frac{\partial Y_{2, m}}{\partial \theta}\quad &&  \nonumber\\
&\qquad \qquad \qquad  -  \left.\left(\frac{\partial^2 \sigma_1}{\partial r^2}+\beta \frac{\partial \sigma_1}{\partial r}\right)\right|_{r=R} Y_{2,m} + \left. \frac{1}{R^2} \frac{\partial \sigma_1}{\partial \theta}\right|_{r= R} \frac{\partial Y_{2, m}}{\partial \theta} \quad && \nonumber \\
&\qquad \qquad \qquad -  \left.\left(\frac{\partial^2 w_1}{\partial r^2}+\beta \frac{\partial w_1}{\partial r}\right)\right|_{r=R} Y_{2,0} + \left. \frac{1}{R^2} \frac{\partial w_1}{\partial \theta}\right|_{r= R} \frac{\partial Y_{2, 0}}{\partial \theta}, \quad && x \in \partial B_R, \\
\label{q_2 equation}
&-\Delta q_2=\mu w_2, \quad && x \in  B_R ,\\
\label{q_2 boundary}
&q_2 = \gamma \frac{1}{R^3}\left( 2 Y_{2,0} Y_{2,m} + Y_{2,0} \Delta_\omega Y_{2,m} + Y_{2,m} \Delta_\omega Y_{2,0} \right) \quad && \nonumber \\ 
&\qquad  -\left. \frac{\partial^2 p_s}{\partial r^2} \right|_{r=R} Y_{2,0} Y_{2,m} -\left. \frac{\partial p_1}{\partial r} \right|_{r=R} Y_{2,m} -\left. \frac{\partial q_1}{\partial r} \right|_{r=R} Y_{2,0}, \quad && x \in \partial B_R.
\end{align}

The solution of problem \eqref{w_1 equation}--\eqref{q_1 boundary} has the following form 
\begin{align*}
w_1 = w_{1,m}(r) Y_{2,m}(\theta,\varphi),   \qquad q_1 = q_{1,m}(r) Y_{2,m}(\theta,\varphi).
\end{align*}
Also, the solution of problem \eqref{sigma_1 equation}--\eqref{p_1 boundary} has the form as follows
\begin{align}
\label{p_1 = p_10 Y_20}
\sigma_1 = \sigma_{1,0}(r) Y_{2,0}(\theta),   \qquad  p_1 = p_{1,0}(r) Y_{2,0}(\theta).
\end{align}
By \eqref{w_1 equation}, \eqref{w_1 boundary}, \eqref{sigma_1 equation} and \eqref{sigma_1 boundary}, both $w_{1,m}$ and $\sigma_{1,0}$ satisfy 
\begin{align*}
& \dfrac{\partial^2 u}{\partial r^2}+\frac{2}{r} \dfrac{\partial u}{\partial r} -\left( 1 + \dfrac{6}{r^2} \right) u = 0, \qquad 0<r<R, \\
&\left. \dfrac{\partial u}{\partial r} \right|_{r=0}=0,\qquad  \left.\left(\dfrac{\partial u}{\partial r}+\beta u\right)\right|_{r=R}=-  \left.\left(\frac{\partial^2 \sigma_s}{\partial r^2}+\beta \frac{\partial \sigma_s}{\partial r}\right)\right|_{r=R}.   
\end{align*}
Hence $w_{1,m}=\sigma_{1,0}$. Similarly, we have $q_{1,m}=p_{1,0}$. Therefore,
\begin{align}
&w_1 = w_{1,m}(r) Y_{2,m}(\theta,\varphi)  = \sigma_{1,0}(r) Y_{2,m}(\theta,\varphi), \nonumber \\ 
\label{q_1 = p_10 Y_2m}
&q_1 = q_{1,m}(r) Y_{2,m}(\theta,\varphi) = p_{1,0}(r) Y_{2,m}(\theta,\varphi).
\end{align}

Substituting \eqref{p_1 = p_10 Y_20} and \eqref{q_1 = p_10 Y_2m} into \eqref{dmu compute1}, \eqref{DRR H Y20 Y2m}, \eqref{w_2 boundary} and \eqref{q_2 boundary} and using \eqref{laplace Ynm}, we obtain 
\begin{align}
\label{DRR H Y20 Y20 1}
&H_{\widetilde{R} \widetilde{R}}\left(0, \mu_2\right)\left[Y_{2,0}, Y_{2, 0}\right] =\left.\frac{\partial^3 p_s}{\partial r^3}\right|_{r=R} Y_{2,0}^2 + 2\left. \frac{\partial^2 p_{1,0}}{\partial r^2}\right|_{r = R} Y_{2, 0}^2 - 2\frac{1}{R^2} \left. p_{1,0} \right|_{r= R} \left( \dfrac{\partial Y_{2, 0}}{\partial \theta} \right)^2  + 2\left. \frac{\partial p_2}{\partial r}\right|_{r = R}, \\
\label{DRR H Y20 Y2m 1}
&H_{\widetilde{R} \widetilde{R}}\left(0, \mu_2\right)\left[Y_{2,0}, Y_{2, m}\right] = \left.\frac{\partial^3 p_s}{\partial r^3}\right|_{r=R} Y_{2,0} Y_{2, m}+ \left. \frac{\partial^2 p_{1,0}}{\partial r^2}\right|_{r = R} Y_{2, 0} Y_{2, m} - \frac{1}{R^2} \left. p_{1,0} \right|_{r= R} \frac{\partial Y_{2, 0}}{\partial \theta} \frac{\partial Y_{2, m}}{\partial \theta} \nonumber \\
&\qquad \qquad \qquad \qquad \qquad \qquad + \left. \frac{\partial^2 p_{1,0}}{\partial r^2}\right|_{r = R} Y_{2, 0} Y_{2, m} - \frac{1}{R^2} \left. p_{1,0} \right|_{r= R} \frac{\partial Y_{2, 0}}{\partial \theta} \frac{\partial Y_{2,m}}{\partial \theta} + \left. \frac{\partial q_2}{\partial r}\right|_{r = R} \nonumber \\
& \qquad \qquad \qquad \qquad 
=\left.\frac{\partial^3 p_s}{\partial r^3}\right|_{r=R} Y_{2,0} Y_{2, m}+ 2\left. \frac{\partial^2 p_{1,0}}{\partial r^2}\right|_{r = R} Y_{2, 0} Y_{2, m} - 2\frac{1}{R^2} \left. p_{1,0} \right|_{r= R} \frac{\partial Y_{2, 0}}{\partial \theta} \frac{\partial Y_{2, m}}{\partial \theta} + \left. \frac{\partial q_2}{\partial r}\right|_{r = R},\\
\label{w_2 boundary 1}
&\left. \left( \frac{\partial w_2}{\partial r}+\beta w_2 \right) \right|_{r=R}= - \left.\left(\frac{\partial^3 \sigma_s}{\partial r^3}+\beta \frac{\partial^2 \sigma_s}{\partial r^2}\right)\right|_{r=R} Y_{2,0} Y_{2,m} + \left. \frac{1}{R^2} \frac{\partial \sigma_s}{\partial r}\right|_{r= R} \frac{\partial Y_{2, 0}}{\partial \theta} \frac{\partial Y_{2, m}}{\partial \theta} \nonumber \\
&\qquad \qquad  \qquad \qquad \qquad  -  \left.\left(\frac{\partial^2 \sigma_{1,0}}{\partial r^2}+\beta \frac{\partial \sigma_{1,0}}{\partial r}\right)\right|_{r=R} Y_{2,0} Y_{2,m} +  \frac{1}{R^2} \left. \sigma_{1,0} \right|_{r= R} \frac{\partial Y_{2, 0}}{\partial \theta} \frac{\partial Y_{2, m}}{\partial \theta} \nonumber \\
&\qquad \qquad \qquad \qquad \qquad -  \left.\left(\frac{\partial^2 \sigma_{1,0}}{\partial r^2}+\beta \frac{\partial \sigma_{1,0}}{\partial r}\right)\right|_{r=R} Y_{2,0} Y_{2,m} +  \frac{1}{R^2} \left. \sigma_{1,0} \right|_{r= R} \frac{\partial Y_{2, 0}}{\partial \theta} \frac{\partial Y_{2, m}}{\partial \theta} \nonumber \\
&\qquad \qquad \qquad \qquad \quad = - \left.\left(\frac{\partial^3 \sigma_s}{\partial r^3}+\beta \frac{\partial^2 \sigma_s}{\partial r^2}\right)\right|_{r=R} Y_{2,0} Y_{2,m} + \left. \frac{1}{R^2} \frac{\partial \sigma_s}{\partial r}\right|_{r= R} \frac{\partial Y_{2, 0}}{\partial \theta} \frac{\partial Y_{2, m}}{\partial \theta} \nonumber \\
&\qquad \qquad  \qquad \qquad \qquad  - 2 \left.\left(\frac{\partial^2 \sigma_{1,0}}{\partial r^2}+\beta \frac{\partial \sigma_{1,0}}{\partial r}\right)\right|_{r=R} Y_{2,0} Y_{2,m} +  2 \frac{1}{R^2} \left. \sigma_{1,0} \right|_{r= R} \frac{\partial Y_{2, 0}}{\partial \theta} \frac{\partial Y_{2, m}}{\partial \theta}, \\
\label{q_2 boundary 1}
&\left. q_2 \right|_{r=R} = \gamma \frac{1}{R^3}\left( 2 Y_{2,0} Y_{2,m} + Y_{2,0} \Delta_\omega Y_{2,m} + Y_{2,m} \Delta_\omega Y_{2,0} \right) \nonumber \\
&\qquad \qquad \quad -\left. \frac{\partial^2 p_s}{\partial r^2} \right|_{r=R} Y_{2,0} Y_{2,m} -\left. \frac{\partial p_{1,0}}{\partial r} \right|_{r=R} Y_{2,0} Y_{2,m} -\left. \frac{\partial p_{1,0}}{\partial r} \right|_{r=R} Y_{2,0} Y_{2,m} \nonumber \\
&\qquad \quad \ = - \gamma \frac{10}{R^3} Y_{2,0} Y_{2,m}-\left. \frac{\partial^2 p_s}{\partial r^2} \right|_{r=R} Y_{2,0} Y_{2,m} - 2\left. \frac{\partial p_{1,0}}{\partial r} \right|_{r=R} Y_{2,0} Y_{2,m}.
\end{align}


The solution of problem \eqref{w_2 equation}, \eqref{w_2 boundary 1}, \eqref{q_2 equation} and \eqref{q_2 boundary 1} has the following form 
\begin{align*}
w_2 = \sum_{n=0}^{4} \sum_{m=-n}^{n} w_{2,n,m}(r) Y_{n,m}(\theta,\varphi),   \qquad q_2 =\sum_{n=0}^{4} \sum_{m=-n}^{n} q_{2,n,m}(r) Y_{n,m}(\theta,\varphi).
\end{align*}
Only $w_{2,1,l}=\left\langle w_2,Y_{1,l} \right\rangle$, $q_{2,1,l}=\left\langle q_2,Y_{1,l} \right\rangle$, $w_{2,2,s}=\left\langle w_2,Y_{2,s} \right\rangle$ and $q_{2,2,s}=\left\langle q_2,Y_{2,s} \right\rangle$ will be used in the later proof. So we just focus on calculating them.

\eqref{w_2 boundary 1} and  \eqref{q_2 boundary 1} imply that $\left. \left(\frac{\partial w_2}{\partial r}+\beta w_2 \right) \right|_{r=R}$ and $\left. q_2 \right|_{r=R}$ are both the linear combination of $Y_{2,0} Y_{2,m}$ and $\frac{\partial Y_{2, 0}}{\partial \theta} \frac{\partial Y_{2, m}}{\partial \theta}$. Using \eqref{product1} and \eqref{product2}, we have 
\begin{align*}
&\left\langle \left.\left(\frac{\partial w_2}{\partial r}+\beta w_2\right)\right|_{r=R},Y_{1, l}  \right\rangle = 0,\qquad \left\langle \left. q_2 \right|_{r=R},Y_{1, l}  \right\rangle = 0,  \qquad  -1\leq l \leq 1,\\
&\left\langle \left.\left(\frac{\partial w_2}{\partial r}+\beta w_2\right)\right|_{r=R},Y_{2, s}  \right\rangle = 0,\qquad \left\langle \left. q_2 \right|_{r=R},Y_{2, s}  \right\rangle = 0, \qquad  s \neq m.
\end{align*}
Together with \eqref{w_2 equation} and \eqref{q_2 equation}, 
we obtain
\begin{align}
\label{q_21l = 0}
&\left\langle w_2, Y_{1,l} \right\rangle = w_{2,1,l}=0, \qquad \left\langle q_2, Y_{1,l} \right\rangle = q_{2,1,l}=0, \qquad -1 \leq l \leq 1, \\
\label{q_22s = 0}
&\left\langle w_2, Y_{2,s} \right\rangle=w_{2,2,s}=0, \qquad \left\langle q_2, Y_{2,s} \right\rangle=q_{2,2,s}=0, \qquad s\neq m.    
\end{align}

Substituting \eqref{p_1 = p_10 Y_20} into \eqref{sigma_2 boundary} and \eqref{p_2 boundary} and using \eqref{laplace Ynm}, we get
\begin{align}
\label{sigma_2 boundary 1}
&\left. \left(\frac{\partial \sigma_2}{\partial r}+\beta \sigma_2 \right) \right|_{r=R} =-\frac{1}{2} \left. \left( \frac{\partial^3 \sigma_s}{\partial r^3} + \beta \frac{\partial^2 \sigma_s}{\partial r^2} \right) \right|_{r=R} Y_{2 ,0}^2  
+ \left.\frac{1}{2 R^2} \frac{\partial \sigma_s}{\partial r}\right|_{r=R}\left(\frac{\partial Y_{2,0}}{\partial \theta}\right)^2 \nonumber \\
&\qquad \qquad \qquad \qquad \qquad - \left.\left(\frac{\partial^2 \sigma_{1,0}}{\partial r^2}+\beta \frac{\partial \sigma_{1,0}}{\partial r}\right)\right|_{r=R} Y_{2 ,0}^2 + \frac{1}{R^2} \left. \sigma_{1,0} \right|_{r=R} \left(\frac{\partial Y_{2,0}}{\partial \theta}\right)^2, \\
\label{p_2 boundary 1}
&\left. p_2 \right|_{r=R} = -\gamma \frac{5}{R^3} Y_{2,0}^2 -  \left.\frac{1}{2} \frac{\partial^2 p_s}{\partial r^2}\right|_{r=R} Y_{2,0}^2 -  \left.\frac{\partial p_{1,0}}{\partial r}\right|_{r=R} Y_{2,0}^2.
\end{align}

Comparing \eqref{w_2 boundary 1} and \eqref{q_2 boundary 1} with \eqref{sigma_2 boundary 1} and \eqref{p_2 boundary 1} and using \eqref{product3}, we derive 
\begin{align}
&\left\langle \left.\left(\frac{\partial w_2}{\partial r}+\beta w_2\right)\right|_{r=R},Y_{2, m}  \right\rangle =2 b_m \left\langle \left.\left(\frac{\partial \sigma_2}{\partial r}+\beta \sigma_2\right)\right|_{r=R},Y_{2, 0}  \right\rangle, \nonumber \\
\label{q_2 Y_2m = 2b_m p_2 Y_20}
&\left\langle \left. q_2 \right|_{r=R},Y_{2, m}  \right\rangle =2 b_m \left\langle \left. p_2 \right|_{r=R},Y_{2, 0}  \right\rangle,
\end{align}
where 
\begin{align} 
\label{b_m}
b_m\left\{\begin{array}{ll}
=1, & m=0, \\
=\dfrac{1}{2}, & m=\pm 1, \\
=-1, & m = \pm 2.
\end{array}\right.
\end{align}
Together with \eqref{w_2 equation}, $w_{2,2,m}$ satisfies 
\begin{align*}
& \dfrac{\partial^2 u}{\partial r^2}+\frac{2}{r} \dfrac{\partial u}{\partial r} -\left( 1 + \dfrac{6}{r^2} \right) u = 0, \qquad 0<r<R, \\
&\left. \dfrac{\partial u}{\partial r} \right|_{r=0}=0,\qquad  \left.\left(\dfrac{\partial u}{\partial r}+\beta u\right)\right|_{r=R} = 2 b_m \left\langle  \left.\left(\dfrac{\partial \sigma_2}{\partial r}+\beta \sigma_2 \right)\right|_{r=R}, Y_{2,0} \right\rangle .   
\end{align*}
On the other hand, \eqref{sigma_2 equation}--\eqref{sigma_2 boundary} imply that $\left\langle \sigma_2,Y_{2, 0}  \right\rangle$ satisfies
\begin{align*}
& \dfrac{\partial^2 u}{\partial r^2}+\frac{2}{r} \dfrac{\partial u}{\partial r} -\left( 1 + \dfrac{6}{r^2} \right) u = 0, \qquad 0<r<R, \\
&\left. \dfrac{\partial u}{\partial r} \right|_{r=0}=0,\qquad  \left.\left(\dfrac{\partial u}{\partial r}+\beta u\right)\right|_{r=R} =  \left\langle  \left.\left(\dfrac{\partial \sigma_2}{\partial r}+\beta \sigma_2 \right)\right|_{r=R}, Y_{2,0} \right\rangle .   
\end{align*}
It follows that $w_{2,2,m} = 2 b_m\left\langle \sigma_2,Y_{2, 0}  \right\rangle$. Similarly, using \eqref{q_2 Y_2m = 2b_m p_2 Y_20}, we have 
\begin{align*}
\left\langle q_2, Y_{2,m} \right\rangle =q_{2,2,m} = 2 b_m\left\langle p_2,Y_{2, 0}  \right\rangle .  
\end{align*}
Together with \eqref{q_21l = 0} and \eqref{q_22s = 0}, we get
\begin{align}
\label{dr q2 Yns}
&\left\langle \left. \frac{\partial q_2}{\partial r}\right|_{r = R} ,Y_{n, s}\right\rangle  = \left\{ \begin{array}{ll}
0,  &n = 1,  \\
0,  &n=2, s \neq m, \\
b_m\left\langle \left. 2\dfrac{\partial p_2}{\partial r}\right|_{r= R} ,Y_{2, 0} \right\rangle,  & n=2, s=m.
\end{array} \right.    
\end{align}


From \eqref{product1}--\eqref{product3}, \eqref{dr3 p_s(R)} and \eqref{p_1 = p_10 Y_20}, we derive
\begin{align}
\label{dr3 ps Yns}
&\left\langle \left.\dfrac{\partial^3 p_s}{\partial r^3}\right|_{r=R} Y_{2,0} Y_{2, m}, Y_{n, s} \right\rangle  = \left\{ \begin{array}{ll}
0,  &n = 1,  \\
0,  &n=2, s \neq m, \\
b_m \left\langle \left.\dfrac{\partial^3 p_s}{\partial r^3}\right|_{r=R} Y_{2,0}^2, Y_{2, 0} \right\rangle,  & n=2, s=m,
\end{array} \right.  \\   
\label{dr2 p1 Yns}
&\left\langle 2\left. \frac{\partial^2 p_{1,0}}{\partial r^2}\right|_{r = R} Y_{2, 0} Y_{2, m},Y_{n, s} \right\rangle  = \left\{ \begin{array}{ll}
0,  &n = 1,  \\
0,  &n=2, s \neq m, \\
b_m \left\langle \left. 2\dfrac{\partial^2 p_{1,0}}{\partial r^2}\right|_{r = R} Y_{2, 0}^2,Y_{2, 0} \right\rangle,  & n=2, s=m,
\end{array} \right.  \\ 
\label{dtheta p1 Yns}
&\left\langle - 2\dfrac{1}{R^2} \left. p_{1,0} \right|_{r= R} \frac{\partial Y_{2, 0}}{\partial \theta} \frac{\partial Y_{2, m}}{\partial \theta} ,Y_{n, s}\right\rangle  = \left\{ \begin{array}{ll}
0,  &n = 1,  \\
0,  &n=2, s \neq m, \\
b_m\left\langle  - 2\dfrac{1}{R^2} \left. p_{1,0} \right|_{r= R} \left( \dfrac{\partial Y_{2, 0}}{\partial \theta} \right)^2,Y_{2, 0} \right\rangle,  & n=2, s=m.
\end{array} \right.  
\end{align}

\eqref{DRR H Y20 Y20 1}--
\eqref{DRR H Y20 Y2m 1} and \eqref{dr q2 Yns}--\eqref{dtheta p1 Yns} imply  
\begin{align}
\label{DRR H Y1 = 0}
&\langle H_{\widetilde{R} \widetilde{R}}(0,\mu_2) \left[ Y_{2,0}, Y_{2, m} \right],Y_{1, l} \rangle = 0, \qquad -1 \leq l \leq 1,\\
\label{DRR H Y2 = 0}
&\langle H_{\widetilde{R} \widetilde{R}}(0,\mu_2) \left[ Y_{2,0}, Y_{2, m} \right],Y_{2, s} \rangle = 0, \qquad s \neq m,\\
\label{DRR H Y2}
&\langle H_{\widetilde{R} \widetilde{R}}(0,\mu_2) \left[ Y_{2,0}, Y_{2, m} \right],Y_{2, m} \rangle  =
 b_m \langle H_{\widetilde{R} \widetilde{R}}(0,\mu_2) \left[ Y_{2,0}, Y_{2, 0} \right],Y_{2, 0} \rangle.
\end{align}

By \eqref{DR H}, we get
\begin{align}
\label{DRmu H Y1 = 0}
&\langle H_{\widetilde{R} \mu}(0,\mu_2) [Y_{2, m}],Y_{1, l} \rangle = 0, \qquad  -1 \leq l \leq 1, \\
\label{DRmu H Y2 = 0}
&\langle H_{\widetilde{R} \mu}(0,\mu_2) [Y_{2, m}],Y_{2, s} \rangle = 0, \qquad  s \neq m, \\
\label{DRmu H}
&\langle H_{\widetilde{R} \mu}(0,\mu_2) [Y_{2, m}],Y_{2, m} \rangle = \langle H_{\widetilde{R} \mu}(0,\mu_2) [Y_{2, 0}],Y_{2, 0} \rangle.  
\end{align}

From \eqref{DRR H Y1 = 0}--\eqref{DRmu H}, \eqref{dmu equation} and \eqref{a},
we have
\begin{align*}
&\langle H_{\widetilde{R} \widetilde{R}}(0,\mu_2) \left[ Y_{2,0}, Y_{2, m} \right],Y_{1, l} \rangle + \langle H_{\widetilde{R} \mu}(0,\mu_2) \mu_2'(0) [Y_{2, m}],Y_{1, l} \rangle = 0, \qquad -1 \leq l \leq 1, \\
&\langle H_{\widetilde{R} \widetilde{R}}(0,\mu_2) \left[ Y_{2,0}, Y_{2, m} \right],Y_{2, s} \rangle + \langle H_{\widetilde{R} \mu}(0,\mu_2) \mu_2'(0) [Y_{2, m}],Y_{2, s} \rangle = 0, \qquad s \neq m, \\
&\langle H_{\widetilde{R} \widetilde{R}}(0,\mu_2) \left[ Y_{2,0}, Y_{2, m} \right],Y_{2, m} \rangle + \langle H_{\widetilde{R} \mu}(0,\mu_2)\mu_2'(0) [Y_{2, m}],Y_{2, m} \rangle \\
&\quad = b_m \langle H_{\widetilde{R} \widetilde{R}}(0,\mu_2) \left[ Y_{2,0}, Y_{2, 0} \right],Y_{2, 0} \rangle + \mu_2'(0) \langle H_{\widetilde{R} \mu}(0,\mu_2) [Y_{2, 0}],Y_{2, 0} \rangle \\
&\quad = -2 b_m \mu_2'(0) \langle H_{\widetilde{R} \mu}(0,\mu_2) [Y_{2, 0}],Y_{2, 0} \rangle + \mu_2'(0) \langle H_{\widetilde{R} \mu}(0,\mu_2) [Y_{2, 0}],Y_{2, 0} \rangle \\
&\quad =-2 \left( b_m - \dfrac{1}{2} \right) \mu_2'(0)  \langle H_{\widetilde{R} \mu}(0,\mu_2) [Y_{2, 0}],Y_{2, 0} \rangle= \left( b_m - \dfrac{1}{2} \right) a.
\end{align*}
Together with \eqref{under Y1m} and \eqref{b_m}, we get the result. 
\end{proof}

\begin{proof}[\textbf{Proof of Proposition \ref{d eigenvalue}}]
By Lemma \ref{8 eigenvalues}, denote $\lambda_i(\varepsilon),i=1,2,\cdots,8$ as the $0$-group eigenvalues of $H_{\widetilde{R}}(\widetilde{R}_2(\varepsilon),\mu_2(\varepsilon))$. Our goal is to compute their derivative.

At first, we use the method in the proof of \cite[VII-\S 1.3, Theorem 1.7]{Kato1976Perturbation} to transform the $0$-group eigenvalues $\lambda_i(\varepsilon),i=1,2,\cdots,8$ of $H_{\widetilde{R}}(\widetilde{R}_2(\varepsilon),\mu_2(\varepsilon))$ to the eigenvalues of $\check{H}(\varepsilon)$ (given by \eqref{check H}) in $Q_0 X^{1+\alpha}$.

Since $0$ is an isolated eigenvalue of $H_{\widetilde{R}}(0,\mu_2)$, the $0$-group eigenvalues of $H_{\widetilde{R}}(\widetilde{R}_2(\varepsilon),\mu_2(\varepsilon))$ are also isolated. We choose a simple closed curve $\Gamma$ to ensure that the $0$-group eigenvalues of $H_{\widetilde{R}}(\widetilde{R}_2(\varepsilon),\mu_2(\varepsilon))$ are all inside $\Gamma$ and other eigenvalues are outside $\Gamma$ for $\varepsilon \in (-\delta,\delta)$. Define $Q(\varepsilon)$ by
\begin{align*}
Q(\varepsilon) = -\frac{1}{2 \pi i} \int_{\Gamma}  \left[ H_{\widetilde{R}}(\widetilde{R}_2(\varepsilon),\mu_2(\varepsilon)) - \zeta I \right]^{-1} d \zeta.
\end{align*}
Similar to the proof of \cite[II-\S 5.4, Theorem 5.4]{Kato1976Perturbation}, we get that $Q(\varepsilon)$ is continuously differentiable to the continuous differentiability of $H_{\widetilde{R}}(\widetilde{R}_2(\varepsilon),\mu_2(\varepsilon))$. By \cite[III-\S 6.4, Theorem 6.17 and its proof]{Kato1976Perturbation}, $Q(\varepsilon)$ is a eigenprojection such that $Q(0)=Q_0$ (given by \eqref{Q_0}) and $X^{1+\alpha}= Q(\varepsilon)X^{1+\alpha} \bigoplus (I-Q(\varepsilon))X^{1+\alpha}$. Then the $0$-group eigenvalues of $H_{\widetilde{R}}(\widetilde{R}_2(\varepsilon),\mu_2(\varepsilon))$ are the eigenvalues of $H_{\widetilde{R}}(\widetilde{R}_2(\varepsilon),$ $\mu_2(\varepsilon)) |_{Q(\varepsilon) X^{1+\alpha}}$.

Although $Q(\varepsilon) X^{1+\alpha}$ is a finite dimensional subspace, it is somewhat inconvenient that this subspace depends on $\varepsilon$. \cite[II-\S 4.2]{Kato1976Perturbation} introduced a continuously differentiable transformation $U(\varepsilon)$ such that $Q(\varepsilon) = U(\varepsilon) Q_0 U(\varepsilon)^{-1}$. Then the subspace $Q(\varepsilon) X^{1+\alpha}$ is transformed into $Q_0 X^{1+\alpha}$ by $U(\varepsilon)$. Define $\check{H}(\varepsilon)$ by
\begin{align}
\label{check H}
\check{H}(\varepsilon) =  U(\varepsilon)^{-1} H_{\widetilde{R}}(\widetilde{R}_2(\varepsilon),\mu_2(\varepsilon)) U(\varepsilon). 
\end{align}
Then the eigenvalue problem for $H_{\widetilde{R}}(\widetilde{R}_2(\varepsilon),$ $\mu_2(\varepsilon)) |_{Q(\varepsilon) X^{1+\alpha}}$ is equivalent to that for $\check{H}(\varepsilon)|_{Q_0 X^{1+\alpha}}$. Hence all the eigenvalue of $\check{H}(\varepsilon)|_{Q_0 X^{1+\alpha}}$ are $\lambda_i(\varepsilon),i=1,2,\cdots,8$.

Next, we compute $\lambda_i'(0)$. Since $\dim Q_0 X^{1+\alpha} = 8$ and $\check{H}(\varepsilon)$ is continuously differentiable in $(-\delta,\delta)$, applying \cite[II-\S 5.4, Theorem 5.4]{Kato1976Perturbation}, we obtain that if $\lambda(\varepsilon)$ is a repeated eigenvalues of $\check{H}(\varepsilon)|_{Q_0 X^{1+\alpha}}$, then $\lambda'(0)$ is a repeated eigenvalues of $Q_0 \check{H}'(0) Q_0$ in $Q_0 X^{1+\alpha}$. Notice that $U(0)=I$. By direct computation, we derive
\begin{align*}
\check{H}'(0) =  \left. \left( U(\varepsilon)^{-1} \right)' \right|_{\varepsilon=0} H_{\widetilde{R}}(0,\mu_2) + \{ H_{\widetilde{R} \widetilde{R}}(0,\mu_2) \left[ Y_{2,0}, \ \cdot \ \right] + H_{\widetilde{R} \mu}(0,\mu_2) \mu_2'(0) \} + H_{\widetilde{R}}(0,\mu_2) U'(0)  
\end{align*}
The fact that $Q_0 H_{\widetilde{R}}(0,\mu_2) = H_{\widetilde{R}}(0,\mu_2)Q_0 =0$ leads to
\begin{align*}
Q_0 \check{H}'(0) Q_0=  Q_0 \{ H_{\widetilde{R} \widetilde{R}}(0,\mu_2) \left[ Y_{2,0}, \ \cdot \ \right] + H_{\widetilde{R} \mu}(0,\mu_2) \mu_2'(0) \} Q_0.    
\end{align*}
Combining Lemma \ref{8 eigenvalues 1}, we get the result.
\end{proof}

\begin{proof}[\textbf{Proof of Theorem \ref{result2}}]
Let $\lambda_1(\varepsilon)$ and $\lambda_2(\varepsilon)$ be the $0$-group eigenvalues of $H_{\widetilde{R}}(\widetilde{R}_2(\varepsilon),\mu_2(\varepsilon))$ given by Proposition \ref{d eigenvalue}. Since $\lambda_1(0)=0$ and $\lambda_2(0)=0$, Proposition \ref{d eigenvalue} implies that for small $|\varepsilon|$, there is at least one $0$-group eigenvalue of $-H_{\widetilde{R}}(\widetilde{R}_2(\varepsilon),\mu_2(\varepsilon))$ with positive value. By \eqref{principle}, we get Theorem \ref{result2}. 
\end{proof}

\section*{Acknowledge}
The research is supported by Guangdong Basic and Applied Basic Research Foundation (Grant No. 2025A1515010922). 


\bibliography{sn-bibliography}

\end{document}